\documentclass{amsart}
\usepackage[latin1]{inputenc}
\usepackage{amssymb}
\usepackage{amsmath}
\usepackage{enumerate}
\usepackage{epsfig}
\usepackage{hyperref}
\usepackage{multicol}
\usepackage{longtable,graphics,multirow,ulem}

\usepackage[all]{xy}

\newtheorem{theorem}{Theorem}[section]
\newtheorem*{theorem*}{Theorem}
\newtheorem{lemma}[theorem]{Lemma}
\newtheorem{proposition}[theorem]{Proposition}
\newtheorem{corollary}[theorem]{Corollary}

\newtheorem{definition}[theorem]{Definition}

\newtheorem{rem}[theorem]{Remark}

\numberwithin{equation}{section}
\newcommand{\rk}{\mbox{rank}}

\newcommand{\Tr}{\mbox{Tr}}

\newcommand{\ra}{\rightarrow}

\newcommand{\eprf}{\hfill $\square$ }

\newcommand{\R}{ \mathbb{R}}
\newcommand{\Z}{\mathbb{Z}}
\newcommand{\Q}{\mathbb{Q}}

\newcommand{\Aut}{\rm{Aut}}
\newcommand{\rt}{\rm{root}}
\newcommand{\tors}{\rm{tors}}

\newcommand{\etalchar}[1]{$^{#1}$}

\makeatletter
\def\blfootnote{\xdef\@thefnmark{}\@footnotetext}

\makeatletter

\newenvironment{figurehere}
{\def\@captype{figure}}
{}
\makeatother

\begin{document}
\title{Classifications of elliptic fibrations of a singular K3 surface}

\author[Bertin]{Marie Jos\'e Bertin}
\email{marie-jose.bertin@imj-prg.fr}

\author[Garbagnati]{Alice Garbagnati}
\email{alice.garbagnati@unimi.it}

\author[Hortsch]{Ruthi Hortsch}
\email{rhortsch@math.mit.edu}

\author[Lecacheux]{Odile Lecacheux}
\email{odile.lecacheux@imj-prg.fr}

\author[Mase]{Makiko Mase}
\email{mtmase@arion.ocn.ne.jp}

\author[Salgado]{Cec\'{i}lia Salgado}
\email{salgado@im.ufrj.br}

\author[Whitcher]{Ursula Whitcher}
\email{whitchua@uwec.edu}

\begin{abstract}
{We classify, up to automorphisms, the elliptic fibrations on the singular K3 surface $X$ whose transcendental lattice is isometric to $\langle 6\rangle\oplus \langle 2\rangle$.}
\end{abstract}

\subjclass[2010]{Primary 14J28, 14J27; Secondary 11G05, 11G42,
14J33} \keywords{K3 surfaces, elliptic fibrations}

\maketitle

\section{Introduction}
We classify elliptic fibrations on the singular K3 surface $X$
associated to the Laurent polynomial
\[
x+\frac{1}{x}+y+\frac{1}{y}+z+\frac{1}{z}+\frac{x}{y}+\frac{y}{x}+\frac{y}%
{z}+\frac{z}{y}+\frac{z}{x}+\frac{x}{z}.%
\]

In order to compute the N\'{e}ron-Severi lattice, the Picard number, and
other basic properties of an algebraic surface, it is useful to identify an elliptic fibration on the surface. Moreover, in view of different applications, one may
be interested in finding all the elliptic fibrations of a certain type. The fibrations of
rank $0$ and maximal torsion lead more easily to the determination of the
$L-$series of the variety \cite{Ber}. Those of positive rank lead to symplectic
automorphisms of infinite order of the variety. Lenstra's Elliptic Curve Method (ECM)
 for finding small factors of large
numbers originally used elliptic curves on $\mathbb{Q}$ with a
torsion-group of order 12 or 16 and rank
$\geq 1$ on $\mathbb{Q}$ \cite{M}, \cite{AM} . One way to obtain infinite families of such curves
is to use fibrations of modular surfaces, as explained by Elkies \cite{El1}.

If the Picard number of a K3 surface is large, there may be an
infinite number of elliptic fibrations, but there is only a finite
number of fibrations up to automorphisms, as proved by Sterk
\cite{St}.  Oguiso used a geometric method to classify elliptic
fibrations in \cite{O}. Some years later, Nishiyama \cite{Nis}
proposed a lattice-theoretic technique to produce such
classifications, recovering Oguiso's results and classifying other
Kummer and
 K3 surfaces. Since then, results of the same type have been obtained by various authors \cite{Ku}, \cite{ES}, \cite{BL}.

Recently, the work of \cite{BKW} described three possible
classifications of elliptic fibrations on a K3 surface, shining a
new light on the meaning of what is a class of equivalence of
elliptic fibrations. In particular, they proposed a
$\mathcal{J}_{1}$-classification of elliptic fibrations up to
automorphisms of the surface and a
$\mathcal{J}_{2}$-classification of the frame lattices of the
fibrations. For our K3 surface, the two classifications coincide.
Thus, it is particularly interesting to exhibit here an
$\mathcal{J}_{2}$-classification by the Kneser-Nishiyama method,
since in general it is not easy to obtain the
$\mathcal{J}_{1}$-classification. This topic will be explained in
detail in Section~\ref{S:classtype}.

Section~\ref{S:presentations} is devoted to a toric presentation
of the surface $X$, following ideas of \cite{KLMSW}, based on the
classification of reflexive polytopes in dimension 3. More
precisely, the Newton polytope of $X$ is in the same class as the
reflexive polytope of index 1529. Since, according to
\cite{KLMSW}, there is an $S_{4}$ action on the vertices of
polytope 1529 and its polar dual, there is a symplectic action  of
$S_{4}$  on $X.$ This action will be described on specific
fibrations. One of them gives the transcendental lattice
$T_{X}=\langle 6 \rangle \oplus\langle 2\rangle.$  We may use
these fibrations to relate $X$ to a modular elliptic surface
analyzed by Beauville in \cite{B}.  We also describe a
presentation of $X$ found in \cite{GS}, which represents $X$ as a
K3 surface with a prescribed abelian symplectic automorphism
group.

 The main results of
the paper are obtained by Nishiyama's method and are summarized in Section~\ref{S:main}, Theorem $4.1.$

\begin{theorem*}
The classification up to automorphisms of the elliptic fibrations on $X$ is given in Table \ref{table: main result}. Each elliptic fibration is given with the Dynkin diagrams characterizing its reducible fibers and the rank and torsion of its Mordell-Weil group.
More precisely, we obtained $52$ elliptic fibrations on $X$, including 17 fibrations of rank $2$ and one of rank $3$.
\end{theorem*}

Due to the high number of different elliptic fibrations, we give
only a few cases of computing the torsion. These cases have been selected to
give an idea of the various methods involved. Notice the case of fibrations \#22
and \#22b exhibiting two elliptic fibrations with the same singular fibers and
torsion but not isomorphic. Corresponding to these different fibrations we
give some particularly interesting Weierstrass models; it is possible to make an exhaustive list.

\textbf {Acknowledgements}

We thank the organizers and all those who supported our project for their
efficiency, their tenacity and expertise.  The authors of the paper have
enjoyed the hospitality of CIRM at Luminy, which helped to initiate a very
fruitful collaboration, gathering from all over the world junior and senior
women, bringing their skill, experience and knowledge  from  geometry and
number theory. Our gratitude goes also to the referee for pertinent remarks and
helpful comments.

A.G is supported by FIRB 2012 ``Moduli Spaces and Their Applications'' and by PRIN 2010--2011 ``Geometria delle variet\`{a} algebriche''.  C.S is supported by FAPERJ (grant E26/112.422/2012).  U.W. thanks the NSF-AWM Travel Grant Program for supporting her visit to CIRM.

\section{Classification of elliptic fibrations on K3 surfaces}\label{S:classtype}
Let $S$ be a smooth complex compact projective surface.
\begin{definition} A surface $S$ is a K3 surface if its canonical bundle  and its irregularity are trivial, that is, if $\mathcal{K}_S\simeq \mathcal{O}_S$ and $h^{1,0}(S)=0$.\end{definition}

\begin{definition}\label{def: elliptic fibration} A flat surjective map $\mathcal{E}:S\ra \mathbb{P}^1$ is called an \emph{elliptic fibration} if:\\ 1) the generic fiber of $\mathcal{E}$ is a smooth curve of genus 1;\\ 2) there exists at least one section $s:\mathbb{P}^1\ra S$ for $\mathcal{E}$.

In particular, we choose one section of $\mathcal{E}$, which we refer to as the zero section. We always denote by $F$ the class of the fiber of an elliptic fibration and by $O$ the curve (and the class of this curve) which is the image of $s$ in $S$.

The group of the sections of an elliptic fibration $\mathcal{E}$ is called the Mordell--Weil group and is denoted by $MW(\mathcal{E})$.\end{definition}

A generic K3 surface does not admit elliptic fibrations, but if the Picard number of the K3 is sufficiently large, it is known that the surface must admit at least one elliptic fibration (see Proposition \ref{prop: e.F. on X with big rho}). On the other hand, it is known that a K3 surface admits a finite number of elliptic fibrations up to automorphisms (see Proposition \ref{prop: finite number of e.f.}).
Thus, a very natural problem is to classify the elliptic fibrations on a given K3 surface. This problem has been discussed in several papers, starting in the Eighties.  There are essentially two different ways to classify elliptic fibrations on K3 surfaces described in  \cite{O} and \cite{Nis}. In some particular cases, a third method can be applied; see \cite{Ku}.
First, however, we must introduce a different problem: ``What does it mean to `classify'  elliptic fibrations?" A deep and interesting discussion of this problem is given in \cite{BKW}, where the authors introduce three different types of classifications of elliptic fibrations and prove that under certain (strong) conditions these three different classifications collapse to a unique one. We observe that it was already known by \cite{O} that in general these three  different classifications do not collapse to a unique one. We now summarize the results by \cite{BKW} and the types of classifications.

\subsection{Types of classifications of elliptic fibrations on K3 surfaces}\label{sec: basic on elliptic fibrations}
In this section we recall some of the main results on elliptic fibrations on K3 surfaces (for example, compare \cite{ScSh}), and we introduce the different classifications of elliptic fibrations discussed in \cite{BKW}.

\subsubsection{The sublattice $U$ and the $\mathcal{J}_0$-classification}\label{sec: U and I0}
Let $S$ be a K3 surface and $\mathcal{E}:S\ra \mathbb{P}^1$ be an
elliptic fibration on $S$. Let $F\in NS(S)$ be the class of the
fiber of $\mathcal{E}$. Then $F$ is a nef divisor which defines
the map $\phi_{|F|}:S\ra \mathbb{P}(H^0(S,F)^*)$ which sends every
point $p\in S$ to $(s_0(p):s_1(p):\ldots :s_r(p))$, where
$\{s_i\}_{i=1,\ldots r}$ is a basis of $H^0(S,F)$, i.e. a basis of
sections of the line bundle associated to the divisor $F$. The map
$\phi_{|F|}$ is the elliptic fibration $\mathcal{E}$. Hence, every
elliptic fibration on a K3 surface is uniquely associated to an
irreducible nef divisor (with trivial self intersection). Since
$\mathcal{E}:S\ra\mathbb{P}^1$ admits a section, there exists a
rational curve which intersects every fiber in one point. Its
class in $NS(S)$ is denoted by $O$ and has the following
intersection properties $O^2=-2$ (since $O$ is a rational curve)
and $FO=1$ (since $O$ is a section). Thus, the elliptic fibration
$\mathcal{E}:S\ra \mathbb{P}^1$ (with a chosen section, as in
Definition \ref{def: elliptic fibration}) is uniquely associated
to a pair of divisors $(F,O)$. This pair of divisors spans a
lattice which is isometric to $U$, represented by the matrix
$\left[\begin{array}{rr}0&1\\1&0\end{array}\right]$, (considering
the basis $F$, $F+O$). Hence each elliptic fibration is associated
to a chosen embedding of $U$ in $NS(S)$.

On the other hand the following result holds:
\begin{proposition}{\rm (\cite[Lemma 2.1]{Ko1} and \cite[Corollary 1.13.15]{Nib})}\label{prop: e.F. on X with big rho}
Let $S$ be a K3 surface, such that there exists a primitive embedding $\varphi:U\hookrightarrow NS(S)$. Then $S$ admits an elliptic fibration.

Let $S$ be a K3 surface with Picard number $\rho(S)\geq 13$. Then, there is a primitive embedding of $U$ in $NS(S)$ and hence $S$ admits at least one elliptic fibration. \end{proposition}



A canonical embedding of $U $ in $NS(S)$ is defined as follows: let us denote by $b_1$ and $b_2$ the unique two primitive vectors of $U$ with trivial self intersection. An embedding of $U$ in $NS(S)$ is called \emph{canonical} if  the image of $b_1$ in $NS(S)$ is a nef divisor and the image of $b_2-b_1$ in $NS(S)$ is an effective irreducible divisor.


The first naive classification of the elliptic fibrations that one
can consider is the classification described above, roughly
speaking: two fibrations are different if they correspond to
different irreducible nef divisors with trivial self
intersections. This essentially coincides with the classification
of the canonical embeddings of $U$ in $NS(S)$.

Following \cite{BKW} we call this classification the $\mathcal{J}_0$-classification of the elliptic fibrations on $S$.

Clearly, it is possible (and indeed likely, if the Picard number
is sufficiently large) that there is an infinite number of
irreducible nef divisors with trivial self intersection and also
infinitely many copies of $U$ canonically embedded in $NS(S)$.
Thus, it is possible that there is an infinite number of
fibrations in curves of genus 1 on $S$ and moreover an infinite
number of elliptic fibrations on $S$.

\subsubsection{Automorphisms and the $\mathcal{J}_1$-classification}\label{subsec: automorphism and J1}
The automorphism group of a variety transforms the variety to itself preserving its structure, but moves points and subvarieties on the variety. Thus, if one is considering a variety with a nontrivial automorphism group, one usually classifies objects on the variety up to automorphisms.

Let $S$ be a K3 surface with a sufficiently large Picard number
(at least 2). Then the automorphism group of $S$ is in general
nontrivial, and it is often of infinite order. More precisely, if
$\rho(S)=2$, then the automorphism group of $S$ is finite if and
only if there is a vector with self intersection either 0 or $-2$
in the N\'eron--Severi group. If $\rho(S)\geq 3$, then the
automorphism group of $S$ is finite if and only if the
N\'eron--Severi group is isometric to a lattice contained in a
known finite list of lattices, cf.\ \cite{Ko2}. Let us assume that
$S$ admits more than one elliptic fibration (up to the
$\mathcal{J}_0$-classification defined above). This means that
there exist at least two elliptic fibrations $\mathcal{E}:S\ra
\mathbb{P}^1$ and $\mathcal{E}':S\ra \mathbb{P}^1$ such that
$F\neq F'\in NS(S)$, where $F$ (resp. $F'$) is the class of the
fiber of the fibration $\mathcal{E}$ (resp. $\mathcal{E}'$). By
the previous observation, it seems very natural to consider
$\mathcal{E}$ and $\mathcal{E}'$ equivalent if there exists an
automorphism of $S$ which sends $\mathcal{E}$ to $\mathcal{E'}$.
This is the idea behind the $\mathcal{J}_1$-classification of the
elliptic fibrations introduced in \cite{BKW}.
\begin{definition}
The $\mathcal{J}_1$-classification of the elliptic fibrations on a K3 surface is the classification of elliptic fibrations up to automorphisms of the surface.
To be more precise: $\mathcal{E}$ is $\mathcal{J}_1$-equivalent to $\mathcal{E}'$ if and only if there exists $g\in \Aut(S)$ such that $\mathcal{E}=\mathcal{E}'\circ g$.
\end{definition}

We observe that if two elliptic fibrations on a K3 surface are equivalent up to automorphism, then all their geometric properties (the type and the number of singular fibers, the properties of the Mordell--Weil group and the intersection properties of the sections) coincide. This is true  essentially by definition, since an automorphism preserves all the ``geometric'' properties of subvarieties on $S$.

The advantages of the $\mathcal{J}_1$-classification with respect
to the $\mathcal{J}_0$-classification are essentially two. The
first is more philosophical: in several contexts, to classify an
object on varieties means to classify the object up to
automorphisms of the variety. The second is more practical and is
based on an important result by Sterk: the
$\mathcal{J}_1$-classification must have a finite number of
classes:

\begin{proposition}{\rm \cite{St}}\label{prop: finite number of e.f.} Up to automorphisms, there exists a finite number of elliptic fibrations on a K3 surface.\end{proposition}

\subsubsection{The frame lattice and the $\mathcal{J}_2$-classification}\label{subsec: frame and J2}
The main problem of the $\mathcal{J}_1$-classification is that it
is difficult to obtain a $\mathcal{J}_1$-classification of
elliptic fibrations on K3 surfaces, since it is in general
difficult to give a complete description of the automorphism group
of a K3 surface and the orbit of divisors under this group. An
intermediate classification can be introduced, the
$\mathcal{J}_2$-classification.  The
$\mathcal{J}_2$-classification is not as fine as the
$\mathcal{J}_1$-classification, and its geometric meaning is not
as clear as the meanings of the classifications introduced above.
However, the $\mathcal{J}_2$-classification can be described in a
very natural way in the context of lattice theory, and there is a
standard method to produce it.

Since the $\mathcal{J}_2$-classification is essentially the
classification of certain lattices strictly related to the
elliptic fibrations, we recall here some definitions and
properties of lattices related to an elliptic fibration.

We have already observed that every elliptic fibration on $S$ is associated to an embedding $\eta:U\hookrightarrow NS(S)$.
\begin{definition} The orthogonal complement of $\eta(U)$ in $NS(S)$, $\eta(U)^{\perp_{NS(S)}}$, is denoted by $W_{\mathcal{E}}$ and called the frame lattice of $\mathcal{E}$.\end{definition}

The frame lattice of $\mathcal{E}$ encodes essentially all the
geometric properties of $\mathcal{E}$, as we explain now. We
recall that the irreducible components of the reducible fibers
which do not meet the zero section generate a root lattice, which
is the direct sum of certain Dynkin diagrams. Let us consider the
root lattice $(W_{\mathcal{E}})_{\rt}$ of $W_{\mathcal{E}}$. Then
the lattice $(W_{\mathcal{E}})_{\rt}$ is exactly the direct sum of
the Dynkin diagram corresponding to the reducible fibers. To be
more precise if the lattice $E_8$ (resp. $E_7$, $E_6$, $D_n$,
$n\geq 4$,  $A_{m}$, $m\geq 3$) is a summand of the lattice
$(W_{\mathcal{E}})_{\rt}$, then the fibration $\mathcal{E}$ admits
a fiber of type $II^*$ (resp. $IV^*$, $III^*$, $I_{n-4}^*$,
$I_{m+1}$). However, the lattices $A_1$ and $A_2$ can be
associated to two different types of reducible fibers, i.e. to
$I_2$ and $III$ and to $I_3$ and $IV$ respectively. We cannot
distinguish between these two different cases using lattice
theory. Moreover, the singular non-reducible fibers of an elliptic
fibration can be either of type $I_1$ or of type $II$.

Given an elliptic fibration $\mathcal{E}$ on a K3 surface $S$, the
lattice $Tr(\mathcal{E}):=U\oplus (W_{\mathcal{E}})_{\rt}$ is
often called the \emph{trivial lattice} (see \cite[Lemma
8.3]{ScSh} for a more detailed discussion).

Let us now consider the Mordell--Weil group of an elliptic
fibration $\mathcal{E}$ on a K3 surface $S$: its properties are
also encoded in the frame $W_{\mathcal{E}}$, indeed
$MW(\mathcal{E})=W_{\mathcal{E}}/(W_{\mathcal{E}})_{\rt}$. In
particular,
$$\rk(MW(\mathcal{E}))=\rk(W_{\mathcal{E}})-\rk((W_{\mathcal{E}})_{\rt})\mbox{
and
}$$
$$(MW(\mathcal{E}))_{\tors}=\overline{(W_{\mathcal{E}})_{\rt}}/(W_{\mathcal{E}})_{\rt},$$
where, for every sublattice $L\subset NS(S)$, $\overline{L}$
denotes the primitive closure of $L$ in $NS(S)$, i.e.\
$\overline{L}:=(L\otimes \Q)\cap NS(S)$.

\begin{definition} The $\mathcal{J}_2$-classification of elliptic fibrations on a K3 surface is the classification of their frame lattices.\end{definition}

It appears now clear that if two elliptic fibrations are
identified by the $\mathcal{J}_2$-classification, they have the
same trivial lattice and the same Mordell--Weil group (since these
objects are uniquely determined by the frame of the elliptic
fibration).

We observe that if $\mathcal{E}$ and $\mathcal{E}'$ are identified
by the $\mathcal{J}_1$-classification, then there exists an
automorphism $g\in Aut(S)$, such that
$\mathcal{E}=\mathcal{E}'\circ g$. The automorphism $g$ induces an
isometry $g^*$ on $NS(S)$ and it is clear that
$g^*:W_{\mathcal{E}}\ra  W_{\mathcal{E}'}$ is an isometry. Thus
the elliptic fibrations $\mathcal{E}$ and $\mathcal{E}'$ have
isometric frame lattices and so are $\mathcal{J}_2$-equivalent.

The $\mathcal{J}_2$-classification is not as fine as the
$\mathcal{J}_1$-classification; indeed, if $h:W_{\mathcal{E}}\ra
W_{\mathcal{E}'}$ is an isometry, a priori there is no reason to
conclude that there exists an automorphism $g\in Aut(S)$ such that
$g^*_{|W_{\mathcal{E}}}=h$; indeed comparing the
$\mathcal{J}_1$-classification given in \cite{O} and the
$\mathcal{J}_2$-classification given in \cite{Nis} for the Kummer
surface of the product of two non-isogenous elliptic curves, one
can check that the first one is more fine than the second one.

The advantage of the $\mathcal{J}_2$-classification sits in its
strong relation with the lattice theory; indeed, there is a method
which allows one to obtain the $\mathcal{J}_2$-classification of
elliptic fibration on several K3 surfaces. This method is
presented in \cite{Nis} and will be described in this paper in
Section \ref{sec: details}.

\subsubsection{Results on the different classification types}
One of the main results of \cite{BKW} is about the relations among
the various types of classifications of elliptic fibrations on K3
surfaces. First we observe that there exists two surjective maps
$\mathcal{J}_0\ra \mathcal{J}_1$ and $\mathcal{J}_0\ra
\mathcal{J}_2$, which are in fact quotient maps (cf.
\cite[Formulae (54) and (57)]{BKW}). This induces a map
$\mathcal{J}_1\ra \mathcal{J}_2$ which is not necessarily a
quotient map.

The \cite[Proposition C']{BKW} gives a bound for the number of different elliptic fibrations up to the $\mathcal{J}_1$-classification, which are identified by the $\mathcal{J}_2$-classification.
As a Corollary the following is proved:
\begin{corollary}\label{cor: J2=J1}{\rm (\cite[Corollary D]{BKW})} Let $S_{(a,b,c)}$ be a K3 surface such that the transcendental lattice of $S$ is isometric to $\left[\begin{array}{ll}2a&b\\b&2c\end{array}\right]$. If $(a,b,c)$ is one of the following $(1,0,1)$, $(1,1,1)$, $(2,0,1)$, $(2,1,1)$, $(3,0,1)$, $(3,1,1)$, $(4,0,1)$, $(5,1,1)$, $(6,1,1)$, $(3,2,1)$, then
$\mathcal{J}_1\simeq\mathcal{J}_2$.\end{corollary}

\subsection{A classification method for elliptic fibrations on K3 surfaces}\label{Nishiyama method in general}

The first paper about the classification of elliptic fibrations on
K3 surfaces is due to Oguiso, \cite{O}. He gives a
$\mathcal{J}_1$-classification of the elliptic fibrations on the
Kummer surface of the product of two non-isogenous elliptic
curves. The method proposed in \cite{O} is very geometric: it is
strictly related to the presence of a certain automorphism (a
non--symplectic involution) on the K3 surface. Since one has to
require that the K3 surface admits this special automorphism, the
method suggested in \cite{O} can be generalized only to certain
special K3 surfaces (see \cite{Kl} and \cite{CG}).

Seven years after the paper \cite{O}, a different method was
proposed by Nishiyama in \cite{Nis}. This method is less geometric
and more related to the lattice structure of the K3 surfaces and
of the elliptic fibrations. Nishiyama applied this method in order
to obtain a $\mathcal{J}_2$-classification of the elliptic
fibrations, both on the K3 surface already considered in \cite{O}
and on other K3 surfaces (cyclic quotients of the product of two
special elliptic curves) to which the method by Oguiso cannot be
applied. Later, in \cite{BL}, the method is used to give a
$\mathcal{J}_2$-classification of elliptic fibrations on a K3
surface whose transcendental lattice is $\langle 4\rangle\oplus
\langle 2\rangle$.

The main idea of Nishiyama's method is the following: we consider a K3 surface $S$ and its transcendental lattice $T_S$. Then we consider a lattice $T$ such that: $T$ is negative definite;  $\rk(T)=\rk(T_S)+4$; the discriminant group and form of $T$ are the same as the ones of $T_S$. We consider primitive embeddings of $\phi:T\hookrightarrow L$, where $L$ is a Niemeier lattice. The orthogonal complement of $\phi(T)$ in $L$ is in fact the frame of an elliptic fibration on $S$.

The classification of the primitive embeddings of $T$ in $L$ for every Niemeier lattice $L$ coincides with the $\mathcal{J}_2$-classification of the elliptic fibrations on $S$. We will give more details on Nishiyama's method in Section \ref{sec: details}.

Since this method is related only to the lattice properties of the surface, a priori one can not expect to find a $\mathcal{J}_1$-classification by using only this method.

Thanks to Corollary \ref{cor: J2=J1}, (see \cite{BKW}) the results
obtained by Nishiyama's method are sometimes stronger than
expected. In particular, we will see that in our case (as in the
case described in \cite{BL}) the classification that we obtain for
the elliptic fibrations on a certain K3 surface using the
Nishiyama's method, is in fact a $\mathcal{J}_1$-classification
(and not only a $\mathcal{J}_2$-classification).

\subsection{Torsion part of the Mordell--Weil group of an elliptic fibration}\label{subset: torsion MW theory}
In Section \ref{sec: explicit computations}, we will classify elliptic fibrations on a certain K3 surface, determining both the trivial lattice and the Mordell--Weil group.  A priori, steps (8) and (9) of the algorithm presented in \ref{sec: details} completely determine the Mordell--Weil group. In any case, we can deduce some information on the torsion part of the Mordell--Weil group by considering only the properties of the reducible fibers of the elliptic fibration. This makes the computation easier, so here we collect some results on the relations between the reducible fibers of a fibration and the torsion part of the Mordell--Weil group.

First, we recall that a section meets every fiber in exactly one
smooth point, so a section meets every reducible fiber in one
point of a component with multiplicity 1 (we recall that the
fibers of type $I_n^*$, $II^*$, $III^*$, $IV^*$ have reducible
components with multiplicity greater than 1). We will call the
component of a reducible fiber which meets the zero section the
\emph{zero component} or \emph{trivial component}.

Every section (being a rational point of an elliptic curve defined over $k(\mathbb{P}^1)$) induces an automorphism of every fiber, in particular of every reducible fiber. Thus, the presence of an $n$-torsion section implies that all the reducible fibers of the fibration admits $\Z/n\Z$ as subgroup of the automorphism group. In particular, this implies the following (well known) result:
\begin{proposition}{\rm (cf. \cite[Section 7.2]{ScSh})}\label{prop: torsion and reducible fibers} Let $\mathcal{E}:S\ra \mathbb{P}^1$ be an elliptic fibration and let $MW(\mathcal{E})_{\tors}$ the torsion part of the Mordell--Weil group.

If there is a fiber of type $II^*$, then $MW(\mathcal{E})_{\tors}=0$.

If there is a fiber of type $III^*$, then $MW(\mathcal{E})_{\tors}\leq (\Z/2\Z)$.

If there is a fiber of type $IV^*$, then $MW(\mathcal{E})_{\tors}\leq (\Z/3\Z)$.

If there is a fiber of type $I_n^*$ and $n$ is an even number, then $MW(\mathcal{E})_{\tors}\leq (\Z/2\Z)^2$.

If there is a fiber of type $I_n^*$ and $n$ is an odd number, then $MW(\mathcal{E})_{\tors}\leq (\Z/4\Z)$.
\end{proposition}
\subsubsection{Covers of universal modular elliptic surfaces}\label{subsec: torsion and modular surfaces}

The theory of universal elliptic surfaces parametrizing elliptic curves with prescribed torsion can also be useful when finding the torsion subgroup of a few elliptic fibrations on the list.  It relies on the following definition/proposition.

\begin{proposition}[{see \cite[2.1.4]{CE} or \cite{ShioModular}}]
Let $\pi: X \rightarrow B$ be an elliptic fibration on a surface $X$. Assume $\pi$ has a section of order $N$, for some $N \in \mathbb{N}$, with $N\geq 4$. Then $X$ is a cover of the universal modular elliptic surface, $\mathcal{E}_N,$ of level $N$.
\end{proposition}
After studying the possible singular fibers of the universal surfaces above, one gets the following.

\begin{proposition}
Let $\mathcal{E}_N$ be the universal modular elliptic surface of level $N $. The following hold:\\
i) If $N\geq 5$ then $\mathcal{E}_N$ admits only semi-stable singular fibers. They are all of type $I_m$ with $m|N$.\\
ii)  The surface $\mathcal{E}_4$ is a rational elliptic surface with singular fibers $I_1^*, I_4, I_1$.\\
\end{proposition}



\subsubsection{Height formula for elliptic fibrations}\label{subset: height}
The group structure of the Mordell--Weil group is the group
structure of the rational points of the elliptic curve defined
over the function field of the basis of the fibration. It is also
possible to equip the Mordell--Weil group of a pairing taking
values in $\Q$, which transforms the Mordell--Weil group to a
$\Q$-lattice. Here we recall the definitions and the main
properties of this pairing. For a more detailed description we
refer to \cite{ScSh} and to the original paper \cite{Shio}.

\begin{definition} Let $\mathcal{E}:S\ra C$ be an elliptic fibration and let $O$ be the zero section. The height pairing is the $\Q$-valued pairing, $<-,->:MW(\mathcal{E})\times MW(\mathcal{E})\ra \Q$ defined on the sections of an elliptic fibration as follows:
$$<P,Q>=\chi(S)+P\cdot O+Q\cdot O-\sum_{c\in \mathcal{C}} contr_c(P,Q),$$

where $\chi(S)$ is the holomorphic characteristic of the surface $S$, $\cdot$ is the intersection form on $NS(S)$, $\mathcal{C}=\{c\in C\mbox{ such that the fiber }\mathcal{E}^{-1}(c)\mbox{ is reducible} \}$ and $contr_c(P,Q)$ is a contribution which depends on the type of the reducible fiber and on the intersection of $P$ and $Q$ with such a fiber as described in \cite[Table 4]{ScSh}.

The value $h(P):=<P,P>=2\chi(S)+2P\cdot O-\sum_{c\in \mathcal{C}}
contr_c(P,P),$ is called the \emph{height} of the section $P$.
\end{definition}

We observe that the height formula is induced by the projection of the intersection form on $NS(S)\otimes \Q$ to the orthogonal complement of the trivial lattice $Tr(\mathcal{E})$ (cf. \cite[Section 11]{ScSh}).
\begin{proposition}{\rm (\cite[Section 11.6]{ScSh})}\label{prop: height and torsion} Let $P\in MW(\mathcal{E})$ be a section of the elliptic fibration $\mathcal{E}:S\ra C$. The section $P$ is a torsion section if and only if $h(P)=0$. \end{proposition}

\section{The K3 surface $X$}\label{S:presentations}

The goal of this paper is the classification of the elliptic fibrations on the unique K3 surface $X$ such that $T_X\simeq \langle 6\rangle\oplus \langle 2\rangle$. This surface is interesting for several reasons, and we will present it from different points of view.

\subsection{A toric hypersurface and the symmetric group $\mathcal{S}_4$}\label{SS:torichypersurface}

Let $N$ be a lattice isomorphic to $\mathbb{Z}^n$.  The dual lattice $M$ of $N$ is given by $\mathrm{Hom}(N, \mathbb{Z})$; it is also isomorphic to $\mathbb{Z}^n$.  We write the pairing of $v \in N$ and $w \in M$ as $\langle v , w \rangle$.

Given a lattice polytope $\diamond$ in $N$, we define its \textit{polar polytope} $\diamond^\circ$ to be $\diamond^\circ = \{w \in $M$ \, | \, \langle v , w \rangle \geq -1 \, \forall \, v \in K\}$.  If $\diamond^\circ$ is also a lattice polytope, we say that $\diamond$ is a reflexive polytope and that $\diamond$ and $\diamond^\circ$ are a mirror pair.  A reflexive polytope must contain $\vec{0}$; furthermore, $\vec{0}$ is the only interior lattice point of the polytope.  Reflexive polytopes have been classified in 1,2,3, and 4 dimensions.  In 3 dimensions, there are 4,319 reflexive polytopes, up to an overall isomorphism preserving lattice structure \cite{KS, KS2}. The database of reflexive polytopes is incorporated in the open-source computer algebra software \cite{Sage}.

Now, consider the one-parameter family of K3 surfaces given by

\begin{equation}\label{E:affineverrillpencil}
x + \frac{1}{x}+ y +\frac{1}{y}+ z +\frac{1}{z}+\frac{x}{y}+\frac{y}{x}+\frac{y}{z}+\frac{z}{y}+\frac{x}{z}+\frac{z}{x}+\lambda.\end{equation}

This family of K3 surfaces was first studied in \cite{Verrill}, where its Picard-Fuchs equation was computed.  A general member of the family has Picard lattice given by $U \oplus \langle 6 \rangle$.

The Newton polytope $\diamond^\circ$ determined by the family of polynomials in Equation~\ref{E:affineverrillpencil} is a reflexive polytope with 12 vertices and 14 facets.  This polytope has the greatest number of facets of any three-dimensional reflexive polytope; furthermore, there is a unique three-dimensional reflexive polytope with this property, up to isomorphism.  In the database of reflexive polytopes found in \cite{Sage}, this polytope has index 1529.

We illustrate $\diamond^\circ$ and its polar polytope $\diamond$ in Figures~\ref{F:2355} and ~\ref{F:1529}.

\begin{multicols}{2}

\begin{figurehere}
\begin{center}
\scalebox{.3}{\includegraphics{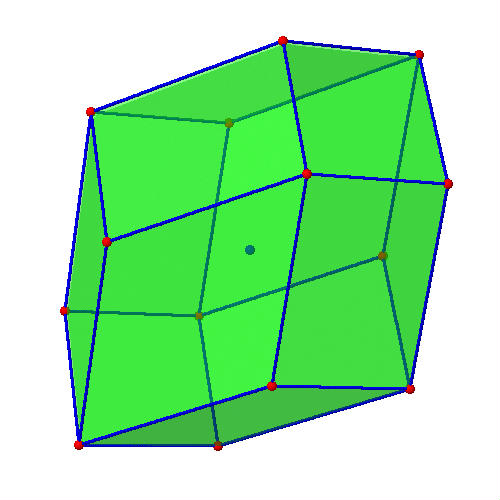}}
\end{center}
\caption{Reflexive polytope 2355}
\label{F:2355}
\end{figurehere}

\begin{figurehere}
\begin{center}
\scalebox{.3}{\includegraphics{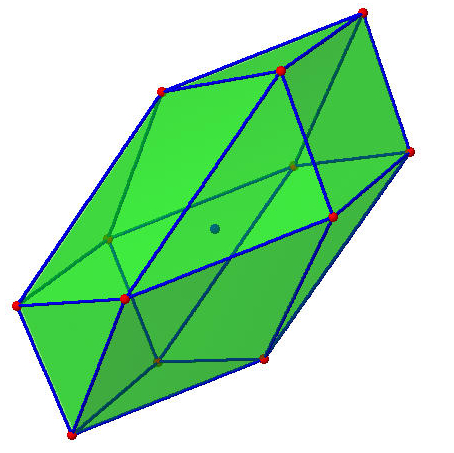}}
\end{center}
\caption{Reflexive polytope 1529}
\label{F:1529}
\end{figurehere}

\end{multicols}

Let us recall some standard constructions and notations involving toric varieties.  A \emph{cone} in $N$ is a subset of the real vector space $N_\mathbb{R} = N \otimes \mathbb{R}$ generated by nonnegative $\mathbb{R}$-linear combinations of a set of vectors $\{v_1, \dots , v_m\} \subset N$.  We assume that cones are strongly convex, that is, they contain no line through the origin.  Note that each face of a cone is a cone.  
A \emph{fan} $\Sigma$ consists of a finite collection of cones such that each face of a cone in the fan is also in the fan, and any pair of cones in the fan intersects in a common face.  We say $\Sigma$ is \emph{simplicial} if the generators of each cone in $\Sigma$ are linearly independent over $\mathbb{R}$.  If every element of $N_\R$ belongs to some cone in $\Sigma$, we say $\Sigma$ is \emph{complete}.  A fan $\Sigma$ defines a toric variety $V_\Sigma$.  If the fan is complete, we may describe $V_\Sigma$ using homogeneous coordinates, in a process analogous to the construction of $\mathbb{P}^n$ as a quotient space of $(\mathbb{C}^*)^n$.  The homogeneous coordinates have one coordinate $z_j$ for each generator of a one-dimensional cone of $\Sigma$.
%

%
We may obtain a fan $R$ from a mirror pair of reflexive polytopes in two equivalent ways.  We may take cones over the faces of $\diamond \subset N_\mathbb{R}$, or we may take the \emph{normal fan} to the polytope $\diamond^\circ \subset M_\mathbb{R}$.  Let $\Sigma$ be a simplicial refinement of $R$ such that the one-dimensional cones of $\Sigma$ are generated by the nonzero lattice points $v_k$, $k = 1 \dots q$, of $\diamond$; we call such a refinement a \emph{maximal projective subdivision}.  Then the variety $V_\Sigma$ is an orbifold.  Then in homogeneous coordinates, we have one coordinate $z_k$ for each nonzero lattice point in $\diamond$.  We may describe the anticanonical hypersurfaces in homogeneous coordinates using polynomials of the form:

\begin{equation}
 p = \sum_{x \in \diamond^\circ \cap M} c_x \prod_{k=1}^q z_k^{\langle v_k, x \rangle + 1}.\end{equation}

\noindent Here the $c_x$ are arbitrary coefficients.  Note that $p$ has one monomial for each lattice point of $\diamond^\circ$.  If the reflexive polytope $\diamond$ is three-dimensional, $V_\Sigma$ is smooth and smooth anticanonical hypersurfaces in $V_\Sigma$ are K3 surfaces (see \cite{CoxKatz} for details).

The orientation-preserving symmetry group of $\diamond$ and $\diamond^\circ$ is the symmetric group $\mathcal{S}_4$.  This group acts transitively on the vertices of $\diamond^\circ$.  As the authors of \cite{KLMSW} observe, by setting the coefficients $c_x$ corresponding to the vertices of $\diamond^\circ$ to 1 and the coefficient corresponding to the origin to a parameter $\lambda$, we obtain a naturally one-parameter family of K3 hypersurfaces with generic Picard rank 19:

\begin{equation}\label{E:verrillpencil}
p = \left( \sum_{x \in \mathrm{vertices}(\diamond^\circ)} \prod_{k=1}^q z_k^{\langle v_k, x \rangle + 1} \right)+\lambda z_1 \dots z_q.\end{equation}

\noindent Equation~\ref{E:verrillpencil} is simply Equation~\ref{E:affineverrillpencil} in homogeneous coordinates.

If we view $\mathcal{S}_4$ as acting on the vertices of $\diamond$ rather than the vertices of $\diamond^\circ$, we obtain a permutation of the homogeneous coordinates $z_k$.  The authors of \cite{KLMSW} show that this action of $\mathcal{S}_4$ restricts to a symplectic action on each K3 surface in the pencil given by Equation~\ref{E:verrillpencil}; in particular, we have a symplectic action of $\mathcal{S}_4$ on $X$.  In the affine coordinates of Equation~\ref{E:affineverrillpencil}, the group action is generated by an element $s_2$ of order 2 which acts by $(x,y,z) \mapsto (1/x, 1/z, 1/y)$ and an element $s_4$ of order 4 which acts by $(x,y,z) \mapsto (x/y, x/z, x)$.

\subsection{The K3 surface $X$}
\begin{definition} Let $X$ be the K3 surface defined by $F=0$, where $F$ is the numerator of
$$x + \frac{1}{x}+ y +\frac{1}{y}+ z +\frac{1}{z}+\frac{x}{y}+\frac{y}{x}+\frac{y}{z}+\frac{z}{y}+\frac{x}{z}+\frac{z}{x}.$$ \end{definition}
The K3 surface $X$ is the special member of the family of K3 surfaces described in \eqref{E:affineverrillpencil} which is obtained by setting $\lambda=0$.\\


We will use three elements of the symplectic group $\mathcal{S}_4$: the three-cycle $s_{3}$ given by $(x,y,z)\mapsto(y,z,x),$ the four-cycle $s_{4}$ and the two-cycle $s_{2}.$

We describe explicitly a first elliptic fibration, which gives the main properties of $X.$

\subsection{A fibration invariant by $s_{3}$}\label{subsec: first explicit fibration}
We use the following factorizations
\begin{align}
F  & =\left(  x+y+z+1\right)  \left(  xy+yz+zx\right)  +\left(
x+y+z-3\right)  xyz\\
F\left(  x+y+z\right)    & =(x+y+z-1)^{2}xyz+(x+y+z+1)(x+y)(y+z)(z+x). \label{FF1}
\end{align}%



If $w=x+y+z,$ we see that $w$ is invariant under the action of $s_{3}.$ If we substitute $w-x-y$ for $z$
 we obtain the equation of an elliptic curve, so the morphism
\begin{align}
 \mathcal{E}:  X &  \rightarrow\mathbb{P}_{w}^{1}\label{u}\\
\left(  x,y,z\right)   &  \mapsto w=x+y+z\nonumber
\end{align}
is an elliptic fibration of $X$.

We use the birational transformation%
\[
x=-\frac{v+(w+1)u}{u(w-3)},\quad y=-\frac{v(w+1)-u^{2}}{v(w-3)}%
\]

with inverse%
\begin{align*}
u&=-((w-3)x+(w+1))((w-3)y+(w+1)), \\
v&=((w-3)x+(w+1))^{2}((w-3)y+(w+1))
\end{align*}

to obtain the Weierstrass equation
\begin{equation}
v^{2}+\left(  w^{2}+3\right)  uv+\left(  w^{2}-1\right)  ^{2}v=u^{3}.
\label{E3}%
\end{equation}

Notice the torsion points $\left(  u=0,v=0\right)$  and $\left(  u=0,v=-(w^{2}%
-1)^{2}\right)  $ of order $3$ and the $3$ points of order $2$ with
$u-$coordinate $\frac{-1}{4}\left(  w^{2}-1\right)  ^{2}$,$-\left(  w-1\right)
^{2}$, and $-\left(  w+1\right)  ^{2}.$

We use also the Weierstrass form%
\begin{equation}
\xi^{2}=\eta\left(  \eta-\left(  w-3\right)  \left(  w+1\right)  ^{3}\right)
\left(  \eta-\left(  w+3\right)  \left(  w-1\right)  ^{3}\right)  \label{We2}%
\end{equation}
with%
\[
u=\frac{1}{4}\left(  \eta-\left(  w^{2}-1\right)  ^{2}\right)  ,\quad
v=\frac{1}{8}\left(  \xi-\left(  w^{2}+3\right)  \eta+\left(  w^{2}-1\right)
^{3}\right).
\]

The singular fibers are of type $I_{6}$ for $w=-1,1,\infty$ and
$I_{2}$ for $w=3,-3,0.$ So the trivial lattice of this fibration
is $\Tr(\mathcal{E})=U\oplus A_5^{\oplus 3}\oplus A_1^{\oplus 3}$.
Hence the Picard number of $X$ is $20$ and
$\rk(MW(\mathcal{E}))=0$. So $X$ is a singular K3 surface.  This
elliptic fibration is contained in the  \cite[Table 2 line 4]{SZ}
and thus its transcendental lattice is $\langle 6\rangle\oplus
\langle 2\rangle$.

Moreover, all the fibers have split mulitiplicative reduction and thus the
N\'{e}ron Severi group is generated by curves defined on $\mathbb{Q}$.

\begin{rem}{\rm We have already observed that $X$ is a special member of the 1-dimensional family of K3 surfaces defined by Equation \ref{E:affineverrillpencil}. Indeed, the transcendental lattice of $X$ is primitively embedded in $U\oplus \langle 6\rangle$ by the vectors $(1,1,0), (0,0,1)$.}\end{rem}

This gives the following proposition:

\begin{proposition}\label{prop: the surface X}
The N\'eron-Severi group of the K3 surface $X$ has rank 20 and is
generated by divisors which are defined over $\mathbb{Q}$. The transcendental
lattice of $X$ is $T_X\simeq \left(
\begin{array}
[c]{cc}%
2 & 0\\
0 & 6
\end{array}
\right).$
\end{proposition}

\begin{rem} {\rm The equation \eqref{We2} is the universal elliptic curve with torsion group
$\mathbb{Z}/2\mathbb{Z}\times\mathbb{Z}/6\mathbb{Z}$ and is in fact equivalent to the equation given in \cite{Kub}.
Thus this fibration can be called modular: we can view the base curve
$(\mathbb{P}_{w}^{1})$ as the modular curve $X_{\Gamma}$ with $\Gamma=\{$
$\left(
\begin{array}
[c]{cc}%
a & b\\
c & d
\end{array}
\right)  \in Sl_{2}\left(  \mathbb{Z}\right)  ,a\equiv1\operatorname{mod}%
6,c\equiv0\operatorname{mod}6,b\equiv0\operatorname{mod}2\}$. 
By \eqref{FF1} we can easily obtain the equation
\[
\left(  w-1\right)  ^{2}xyz+\left(  w+1\right)  \left(  x+y\right)  \left(
y+z\right)  \left(  z+x\right)  =0
\]
and realize $X$ by a base change of the modular rational elliptic
surface $\mathcal{E}_6$ described by Beauville in \cite{B}. We can
prove that on the fiber, the automorphism $s_{3}$ corresponds to
adding a $3$-torsion point.}
\end{rem}

\subsection{A fibration invariant by $s_{4}$}\label{subsec: second explicit fibration}

If $t=\frac{y}{zx},$ we see that $t$ is invariant under the action
of $s_{4}.$ Substituting $tzx$ for $y$ in $F$, we obtain the
equation of an elliptic curve. Using standard transformations (as
in \cite{Ku}(39.2), \cite{AM} or \cite{C}) we obtain the
Weierstrass model
\[
v^{2}=u\left(  u^{2}-2t\left(  t^{2}+1\right)
u+t^{2}(t+1)^{4}\right).
\]

The point $Q_{t}=\left(  u=t(t+1)^{2},v=2t^{2}(t+1)^{2}\right)  $ is of order
$4.$

The point $P_{t}=\left(  u=(t+1\right)  ^{2},v=(t^{2}+1)\left(  t+1\right)
^{2})$ is of infinite order.

The singular fibers are $2I_{1}^{\ast}\left(
t=0,\infty\right)+I_{8}\left(  t=-1\right)+2 I_{1}\left(
t^{2}+t+1=0\right).$

One can prove that on the fiber, $s_{4}$ corresponds to the translation by a $4-$torsion
point. Moreover, the translation by the point $P_{t}$ defines an automorphism of infinite
order on $X$.
\begin{rem}{\rm
If we compute the height of $P_{t}$ we can show, using Shioda formula, \cite{Shio} that
$P_{t}$ and  $Q_{t}$ generate the Mordell-Weil group. }
\end{rem}

\subsection{A fibration invariant by $s_{2}$}\label{subsec: third explicit fibration}

If $r=\frac{y}{z}$ , we see that $r$ is invariant under the action of $s_{2}.$
Substituting $rz$ for $y$ we obtain the equation of an elliptic curve and the following Weierstrass model
\[
v^{2}-\left(  r^{2}-1\right)  vu=u\left(  u-2r\left(  r+1\right)  \right)
\left(  u-2r^{2}\left(  r+1\right)  \right).
\]
\noindent
The point $\left(  0,0\right)  $ is a two-torsion point.
The point $\left(  2r\left(  r+1\right)  ,0\right)  $ is of infinite order.

The singular fibers are

$2I_{6}\left(  r=0,\infty\right)+I_{0}^{\ast}\left(  r=-1\right)+
I_{4}\left(  r=1\right)+2I_{1}\left(  r^{2}-14r+1=0\right)$.

\begin{rem}\label{rem: elliptic fibration Elkies}{\rm From Elkies results \cite{El} \cite{Schu} there is a unique K3 surface
$X/\mathbb{Q}$ with N\'{e}ron-Severi group of rank $20$ and discriminant $-12$ that consists entirely of classes of divisors defined over $\mathbb{Q}$. Indeed it is $X$. Moreover, in \cite{El} a Weierstrass equation for an elliptic fibration on $X$ defined over $\Q$, is given:
\[
y^{2}=x^{3}-75x+(4t-242+\frac{4}{t}).
\]}
\end{rem}

\begin{rem}{\rm The surface $X$ is considered also in a slightly different context in \cite{GS} because of its relation with the study of the moduli space of K3 surfaces with
a symplectic action of a finite abelian group. Indeed, the aim of
the paper \cite{GS} is to study elliptic fibrations
$\mathcal{E}_G:S_G\ra \mathbb{P}^1$ on K3 surfaces $S_G$ such that
$MW(\mathcal{E}_G)=G$ is a torsion group. Since the translation by
a section is a symplectic automorphism of $S_G$, if
$MW(\mathcal{E}_G)=G$, then $G$ is a group which acts
symplectically on $S_G$. In \cite{GS} it is shown how one can
describe both a basis for the N\'eron--Severi group of $S_G$ and
the action induced by the symplectic action of $G$ on this basis.
In particular, one can directly compute the lattices $NS(S_G)^G$
and $\Omega_G:=NS(S_G)^{\perp}$. The latter does not depend on
$S_G$ but only on $G$ and its computation plays a central role in
the description of the moduli space of the K3 surfaces admitting a
symplectic action of $G$ (see \cite{Nisym} and \cite{GS}).  In
particular, the case $G=\Z/6\Z\times Z/2\Z$ is considered: in this
case the K3 surface $S_G$ is $X$, and the elliptic fibration
$\mathcal{E}_G$ is \eqref{u}. Comparing the symplectic action of
$\Z/6\Z\times\Z/2\Z$ on $X$ given in \cite{GS} with the symplectic
group action of $\mathcal{S}_4$ described in
\S~\ref{SS:torichypersurface}, we find that the two groups
intersect in the subgroup of order 3 generated by the map $s_3$
given by $(x,y,z) \mapsto (z,y,x)$. }\end{rem}

\section{Main result}\label{S:main}
This section is devoted to the proof of our main result:

\begin{theorem} The classification up to automorphisms of the elliptic fibrations on $X$ is given in Table \ref{table: main result}. Each elliptic fibration is given with the Dynkin diagrams characterizing its reducible fibers and the rank and torsion of its Mordell-Weil group.
More precisely, we obtained $52$ elliptic fibrations on $X$, including 17 fibrations of rank $2$ and one of rank $3$.
\end{theorem}

We denote by $r$ the $\rk(MW(\mathcal{E}))$, and we use Bourbaki notations for $A_n,D_n,E_k$ as in \cite{BL}.

\begin{longtable}{|c|c|c|c|c|c|c|c|c|c|}\caption{The elliptic fibrations of $X$}\label{table: main result} \\
\hline
$ L_{\text{root}}$ & No.  &  &  &$N_{\text{root}}$ & $r$ & $MW(\mathcal{E})_{\tors}$ \\ \hline
  $E_8^3$  \\ \hline

&$1$   & $A_5\oplus A_1 \subset E_8$ &   & $A_1 E_8  E_8$ & $1$ &
$(0)$\\ \hline &$2$   &  $A_1\subset E_8$ & $A_5 \subset E_8$  &
$A_1 A_2 E_7  E_8$ & $0$ & $(0)$\\ \hline
  $E_8 D_{16}$    \\ \hline

& $3$    & $A_5\oplus A_1 \subset E_8$ &   & $A_1 D_{16}$ & $1$ &
$\mathbb Z /{2 \mathbb Z}$\\ \hline & $4$    & $A_5\oplus A_1
\subset D_{16}$ &   & $ A_1 D_{8}E_8$ & $1$ & $(0)$\\ \hline & $5$
& $A_5\subset E_8$ & $A_1 \subset D_{16}$  & $ A_1A_1 A_2 D_{14}$
& $0$ & $\mathbb Z /{2 \mathbb Z}$\\ \hline &$6$     & $A_1\subset
E_8$ & $A_5 \subset D_{16}$  & $E_7 D_{10}$ & $1$ & $(0)$\\ \hline
  $E_7^2 D_{10}$ \\ \hline
  &$7$  &    $A_5\oplus A_1 \subset E_7$ &   & $E_7 D_{10}$ & $1$ & $\mathbb Z /{2 \mathbb Z}$\\ \hline
  &$8$ &   $A_5\oplus A_1 \subset D_{10}$ &   & $A_1A_1A_1E_7 E_7$ & $1$ & $\mathbb Z /{2 \mathbb Z}$\\ \hline
  &$9$ &    $A_1\subset E_7$ & $A_5 \subset E_7$  & $ D_6 A_1 D_{10}$ &$1$ & $\mathbb Z/{2 \mathbb Z}$\\ \hline
  &$10$&    $A_1\subset E_7$ & $A_5 \subset E_7$  & $ D_6 A_2 D_{10}$ &$0$ &  $\mathbb Z /{2 \mathbb Z}$\\ \hline

   &$11$ &    $A_5\subset E_7$ & $A_1 \subset D_{10}$  & $A_1 A_1 D_8 E_7$ & $1$ & $\mathbb Z /{2 \mathbb Z}$\\ \hline

   &$12$ &    $A_5\subset E_7$ & $A_1 \subset D_{10}$  & $A_1 A_2 D_8 E_7$ & $0$ & $\mathbb Z /{2 \mathbb Z}$\\ \hline

  &$13$ &    $A_1\subset E_7$ & $A_5 \subset D_{10}$  & $E_7 D_6 D_4$ & $1$ & $\mathbb Z /{2 \mathbb Z}$\\ \hline

 $E_7 A_{17}$     \\ \hline
  &$14$ &    $A_5\oplus A_1 \subset E_7$ &   & $A_{17}$ & $1$ & $\mathbb Z /{3 \mathbb Z}$\\ \hline
  &$15$ &    $A_5\oplus A_1 \subset A_{17}$ &   & $A_9E_7$ & $2$ & $(0)$\\ \hline
  &$16$ &    $A_5\subset E_7$ & $A_1 \subset A_{17}$  & $A_1 A_{15}$ & $2$ & $(0)$\\ \hline
 &$17$ &    $A_5\subset E_7$ & $A_1 \subset A_{17}$  & $A_2 A_{15}$ & $1$ & $(0)$\\ \hline
 &$18$ &    $A_1\subset E_7$ & $A_5 \subset A_{17}$  & $D_6 A_{11}$ & $1$ & $(0)$\\ \hline
 $D_{24}$  \\ \hline
&$19$ &    $A_5\oplus A_1 \subset D_{24}$ &   & $A_1 D_{16}$ & $1$
& $(0)$\\ \hline
 $D_{12}^2$   \\ \hline
&  $20$ &   $A_5\oplus A_1 \subset D_{12}$ &   & $A_1 D_4 D_{12}$
& $1$ & $\mathbb Z /{2 \mathbb Z}$\\ \hline

&$21$ &    $A_1 \subset D_{12}$ & $A_5 \subset D_{12}$  & $A_1 D_{10} D_{6}$ & $1$ & $\mathbb Z /{2 \mathbb Z}$\\ \hline
  $D_8^3$     \\ \hline
 &$22$ &   $A_5\oplus A_1 \subset D_{8}$ &   & $A_1 D_8 D_{8}$ & $1$ & $\mathbb Z /{2 \mathbb Z}$\\ \hline
 &$22(b)$   & $A_5\oplus A_1 \subset D_{8}$ &   & $A_1 D_8 D_{8}$ & $1$ & $\mathbb Z /{2 \mathbb Z}$\\ \hline
 & $23$   & $A_1 \subset D_{8}$ & $A_5 \subset D_{8}$  & $A_1^3 D_{6} D_{8}$ & $1$ & $(\mathbb Z /{2\mathbb{Z}})^2 $\\ \hline
 $D_9 A_{15}$ \\ \hline
&  $24$    & $A_5\oplus A_1 \subset D_{9}$ &   & $ A_1 A_{15}$ &
$2$ & $\mathbb Z /{2 \mathbb Z}$\\ \hline & $25$    & $A_5\oplus
A_1 \subset A_{15}$ &   & $ D_9 A_{7}$ & $2$ & $(0)$\\ \hline &
$26$   & $A_5 \subset D_{9}$ & $A_1 \subset A_{15}$  & $A_3
A_{13}$ & $2$ & $(0)$\\ \hline & $27$   & $A_1 \subset D_{9}$ &
$A_5 \subset A_{15}$  & $A_1A_9 D_7$ & $1$ & $(0)$\\ \hline
$E_6^4$  \\ \hline &$28$   & $A_1 \subset E_6$ & $A_5 \subset E_6$
& $A_1A_5 E_6^2$ & $0$ & $\mathbb Z/3\mathbb Z$\\ \hline
 $ A_{11} E_6 D_7$  \\ \hline
&$29$   & $A_5\oplus A_1 \subset A_{11}$ &   & $A_3D_7E_6$ & $2$ &
$(0)$\\ \hline & $30$   & $A_5 \subset E_6$ & $A_1 \subset D_7$  &
$A_1^2A_{11}D_5$ & $0$ & $\mathbb Z/4\mathbb Z$\\ \hline &$31$   &
$A_5 \subset E_6$ & $A_1 \subset A_{11}$  & $A_1A_9D_7 $ & $1$ &
$(0)$\\ \hline &$32$   & $A_1 \subset E_6$ & $A_5 \subset D_7$  &
$ A_5A_{11}$ & $2$ & $\mathbb Z/3\mathbb Z$\\ \hline &$33$   &
$A_5 \subset D_7$ & $A_1 \subset A_{11}$  & $ A_9E_6$ & $3$ &
$(0)$\\ \hline &$34$   & $A_1 \subset D_7$ & $A_5 \subset A_{11}$
& $A_1 A_5D_5E_6$ & $1$ & $(0)$\\ \hline &$35$   & $A_1 \subset
E_6$ & $A_5 \subset A_{11}$  & $ A_5^2D_7$ & $1$ & $(0)$\\ \hline

 $D_6^4$     \\ \hline
&$36$   & $A_1 \subset D_6$ & $A_5 \subset D_6$  & $ A_1 D_4 D_6^2$ & $1$ & $(\mathbb Z /{2\mathbb{Z}})^2$\\ \hline
 $D_6 A_9^2$    &\\ \hline
&$37$    &  $A_5\oplus A_1 \subset A_9$&  & $ A_1 A_9D_6$ & $2$ &
$\mathbb Z/2\mathbb Z$\\ \hline &$38$ & $A_5 \subset A_9$ & $A_1
\subset A_9$ & $A_3A_7D_6$ & $2$ & $(0)$\\\hline

&$39$ & $A_5 \subset A_9$ & $A_1 \subset D_6$ & $A_1A_3A_9D_4$ & $1$ & $\mathbb Z/2\mathbb Z$\\\hline
&$40$ & $A_1 \subset A_9$ & $A_5 \subset D_6$ & $A_7A_9$ & $2$ & $(0)$\\\hline

$D_5^2 A_7^2$     \\ \hline
&$41$ & $A_5\oplus A_1\subset A_7$ && $A_7D_5^2$ & $1$ & $\mathbb Z/4\mathbb Z$\\ \hline
&$42$ & $A_5 \subset A_7$ & $A_1 \subset A_7$ & $A_1A_5D_5^2$ & $2$ & $(0)$ \\ \hline
&$43$  & $A_5 \subset A_7$ & $A_1 \subset D_5$ & $A_1^2A_3A_7D_5$ & $1$ & $\mathbb Z/4\mathbb Z$ \\ \hline

$A_8^3$  \\ \hline &$44$ & $A_5 \oplus A_1 \subset A_8$ & &
$A_8^2$ & $2$ & $\mathbb Z /3 \mathbb Z$\\ \hline &$45$  & $A_1
\subset A_8$ & $A_5 \subset A_8$ & $A_2A_6A_8$ & $2$ & $(0)$\\
\hline

$A_{24}$ \\ \hline
&$46$ & $A_5 \oplus A_1 \subset A_{24}$ & & $A_{16}$ & $2$ & $(0)$\\ \hline

$A_{12}^2$  \\ \hline
&$47$  & $A_5 \oplus A_1 \subset A_{12}$ & & $A_4A_{12}$ & $2$ & $(0)$\\ \hline
&$48$  & $A_5 \subset A_{12}$ & $A_1 \subset A_{12}$ & $ A_6A_{10}$ & $2$ & $(0)$\\ \hline

$D_4A_5^4$ \\ \hline
&$49$  & $A_5=A_5$ & $A_1 \subset A_5$ & $A_3A_5^2D_4$ & $1$ & $\mathbb Z/2\mathbb Z$\\ \hline
&$50$  & $A_5=A_5$ & $A_1 \subset D_4$ & $A_1^3A_5^3$ & $0$ & $\mathbb Z/2\mathbb Z\times \mathbb Z/6 \mathbb Z$\\ \hline

$A_6^4$  \\ \hline
&$51$ & $A_5 \subset A_6$ & $A_1 \subset A_6$ & $A_4A_6^2$ & $2$ & $(0)$\\ \hline
\end{longtable}

{\it Outline of the proof.\ } The proof consists of an application
of Nishiyama's method: the details of this method will be
described in Section \ref{sec: details}. Its application to our
case is given in Section \ref{sec: explicit computations}. The
application of Nishiyama's method gives us a
$\mathcal{J}_2$-classification, which coincides in our cases with
a classification up to automorphisms of the surface by Corollary
\ref{cor: J2=J1}.



\begin{rem}{\rm The fibration given in Section \ref{subsec: first explicit fibration} is \# 50 in Table~\ref{table: main result},  the one given in Section \ref{subsec: second explicit fibration} is \# 41, the one given in Section \ref{subsec: third explicit fibration} is \# 49, the one given in Remark \ref{rem: elliptic fibration Elkies} is \#1. The fibrations of rank 0 may be found also in \cite{SZ}.}\end{rem}
\begin{rem}{\rm We observe that there exists a primitive embedding of $T_X\simeq \langle 6\rangle \oplus\langle 2\rangle$ in $U(2)\oplus\langle 2\rangle$ given by the vectors $\langle(1,1,1),(0,-1,1)\rangle$. Thus, $X$ is a special member of the 1-dimensional family of K3 surfaces whose transcendental lattice is isometric to $U(2)\oplus\langle 2\rangle$. The elliptic fibrations on the generic member $Y$ of this family have already been classified (cf. \cite{CG}), and indeed the elliptic fibrations in Table \ref{table: main result} specialize the ones in \cite[Table 4.5 and Section 8.1 case $r=19$]{CG}, either because the rank of the Mordell--Weil group increases by 1 or because two singular fibers glue together producing a different type of reducible fiber.}\end{rem}

\subsection{Nishiyama's method in detail: an algorithm}\label{sec: details}
This section is devoted to a precise description of Nishiyama's method. Since the method is very well described both in the original paper \cite{Nis} and in some other papers where it is applied,  e.g. \cite{BL} and \cite{BKW}, we summarize it in an algorithm which allows us to compute all the results given in Table \ref{table: main result}. In the next section we will describe in detail some peculiar cases, in order to show how this algorithm can be applied.

\begin{definition} A Niemeier lattice is an even unimodular negative definite lattice of rank 24. \end{definition}
There are 24 Niemeier lattices. We will denote by $L$ an arbitrary Niemeier lattice. Each of them corresponds uniquely to its root lattice $L_{\rt}$.

In Table \ref{table: generators of Niemeier lattice} we list the Niemeier lattices, giving both the root lattices of each one and a set of generators for $L/L_{\rt}$. To do this we recall the following notation, introduced in \cite{BL}:
\begin{eqnarray*}
\begin{array}{|ll||ll|}\hline \alpha_n&=\frac {1}{n+1} \sum_{j=1}^{n}(n-j+1)a_j&
\delta_{l} & = \frac {1}{2}\left( \sum_{i=1}^{l-2}id_i+\frac {1}{2}(l-2)d_{l-1}+\frac {1}{2}ld_l \right)\\
\overline{\delta}_l & =\sum_{i=1}^{l-2} d_i+\frac {1}{2}(d_{l-1}+d_l)&
\tilde{\delta}_{l} & =\frac {1}{2}\left( \sum_{i=1}^{l-2}id_i+\frac {1}{2}ld_{l-1}+\frac {1}{2}(l-2)d_l  \right)\\
\eta_6 &=-\frac {2e_1+3e_2+4e_3+6e_4+5e_5+4e_6}{3}&
\eta_7 &=-\frac {(2e_1+3e_2+4e_3+6e_4+5e_5+4e_6+3e_7)}{2}\\ \hline
\end{array}\end{eqnarray*}

\begin{longtable}{|c|l|}
\caption{The Niemeier lattices $L$:  $L_{\rt}$ and $L/L_{\rt}$}\label{table: generators of Niemeier lattice}\\
\hline
$L_{\rt}$&$L/L{\rt}$\\
\hline
 $E_8^{ 3}$    &  $\langle(0) \rangle$\\ \hline
 $E_8 D_{16}$    &  $\Z/2\Z=\langle \delta_{16} \rangle $\\ \hline
 $E_7^{ 2} D_{10}$    &  $(\Z/2\Z)^2=\langle \eta_7^{(1)}+\delta_{10}, \eta_7^{(1)}+\eta_7^{(2)}+\overline{\delta_{10}} \rangle $ \\ \hline
 $E_7 A_{17}$    &  $\Z/6\Z=\langle \eta_7+3\alpha_{17} \rangle$\\ \hline
 $D_{24}$    &  $\Z/2\Z=\langle \delta_{24} \rangle $ \\ \hline
 $D_{12}^{ 2}$    &  $(\Z/2\Z)^2=\langle \delta_{12}^{(1)}+\overline{\delta_{12}^{(2)}},\overline{\delta_{12}^{(1)}}+\delta_{12}^{(2)} \rangle$ \\ \hline
 $D_8^{ 3}$    &  $(\Z/2\Z)^3=\langle \delta_8^{(1)}+\overline{\delta}_8^{(2)}+\overline{\delta}_8^{(3)}, \overline{\delta}_8^{(1)}+\delta_8^{(2)}+\overline{\delta}_8^{(3)},\overline{\delta}_8^{(1)}+\overline{\delta}_8^{(2)}+\delta_8^{(3)}\rangle $ \\ \hline
  $D_9 A_{15}$    &  $\Z/8\Z=\langle \delta_9+2\alpha_{15} \rangle$ \\ \hline
  $E_6^{ 4}$    &  $(\Z/3\Z)^2=\langle \eta_6^{(1)}+\eta_6^{(2)}+\eta_6^{(3)}, -\eta_6^{(1)}+\eta_6^{(3)}+\eta_6^{(4)} \rangle$\\ \hline
   $ A_{11} E_6 D_7$    &  $ \Z/12\Z=\langle \alpha_{11}+\eta_6+\delta_7 \rangle$ \\ \hline
 $D_6^{ 4}$    &  $\begin{array}{ll}(\Z/2\Z)^4=&\langle {\delta}_6^{(2)}+\overline{\delta}_6^{(3)}+\tilde{\delta}_6^{(4)},\overline{\delta_6^{(1)}}+\tilde{\delta}_6^{(2)}+{\delta}_6^{(4)},\\&\delta_6^{(1)}+\overline{\delta}_6^{(2)}+\tilde{\delta}_6^{(4)},\delta_6^{(1)}+\tilde{\delta}_6^{(3)}+\overline{\delta}_6^{(4)}\rangle\end{array}$\\ \hline
 $D_6 A_9^{ 2}$    &  $\Z/2\Z\times \Z/10\Z=\langle \widetilde{\delta_6}+5\alpha_9^{(2)}, \delta_6+\alpha_9^{(1)}+2\alpha_9^{(2)} \rangle$ \\ \hline
 $D_5^{ 2} A_7^{ 2}$    &  $\Z/4\Z\times\Z/8\Z=\langle \delta_5^{(1)}+\delta_5^{(2)}+2\alpha_7^{(1)},\delta_5^{(1)}+2\delta_5^{(2)}+\alpha_7^{(1)}+\alpha_7^{(2)} \rangle$ \\ \hline
 $ A_8^{ 3}$    &  $\Z/3\Z\times \Z/9\Z=\langle 3\alpha_8^{(1)}+3\alpha_8^{(2)},\alpha_8^{(1)}+2\alpha_8^{(2)}+2\alpha_8^{(3)} \rangle$
  \\ \hline
 $ A_5^{ 4}D_4$    &  $\begin{array}{ll}(\Z/6\Z)^2\times \Z/2\Z=&\langle 5\alpha_5^{(1)}+2\alpha_5^{(2)}+\alpha_5^{(3)}+\overline{\delta}_4,\\& 5\alpha_5^{(1)}+3\alpha_5^{(2)}+2\alpha_5^{(3)}+4\alpha_5^{(4)}+\delta_4,\\&3\alpha_5^{(1)}+3\alpha_5^{(4)}+\tilde{\delta}_4 \rangle\end{array}$\\ \hline
  $ A_6^{ 4}$    &  $(\Z/7\Z)^2=\langle \alpha_6^{(1)}+2\alpha_6^{(2)}+\alpha_6^{(3)}+6\alpha_6^{(4)},\alpha_6^{(1)}+6\alpha_6^{(2)}+2\alpha_6^{(3)}+\alpha_6^{(4)} \rangle$ \\ \hline
 $ D_{4}^{ 6}$    & $(\Z/2\Z)^6=\langle \delta_4^{(1)}+\delta_4^{(i)},\ i=1,\ldots 5,\ \sum_{i=1}^6\widetilde{\delta_4^{(i)}} \rangle$ \\ \hline
 $ A_{24}$    & $\Z/5\Z=\langle5\alpha_{24}\rangle$ \\ \hline
 $ A_{12}^{ 2}$    & $\Z/13\Z=\langle 2\alpha_{13}^{(1)}+3\alpha_{13}^{(2)}\rangle$ \\ \hline
$ A_4^{ 6}$    & $\begin{array}{ll}(\Z/5\Z)^3=&\langle
\alpha_4^{(1)}+\alpha_4^{(2)}+\alpha_4^{(3)}+4\alpha_4^{(4)}+4\alpha_4^{(5)},\\&
\alpha_4^{(1)}+\alpha_4^{(2)}+4\alpha_4^{(3)}+\alpha_4^{(5)}+4\alpha_4^{(6)},\\
&\alpha_4^{(1)}+4\alpha_4^{(3)}+\alpha_4^{(4)}+4\alpha_4^{(5)}+\alpha_4^{(6)}\rangle\end{array}$ \\ \hline

$ A_3^{ 8}$    & $\begin{array}{ll}(\Z/4\Z)^4=&\langle 3\alpha_3^{(1)}+\sum_{i=2}^8c_i\alpha_3^{(i)}\mbox{ such that }(c_2,\ldots c_8)\\&\mbox{ is a cyclic permutation of }(2001011)\rangle\end{array}$ \\ \hline
$ A_2^{ 12}$    & $\begin{array}{ll}(\Z/3\Z)^6=&\langle 2\alpha_2^{(1)}+\sum_{i=2}^{12}c_i\alpha_2^{(i)}\mbox{ such that } (c_2,\ldots c_{12})\\&\mbox{ is a cyclic permutation of }(11211122212)\rangle\end{array}$ \\ \hline
$ A_1^{ 24}$    & $\begin{array}{ll}(\Z/2\Z)^{12}=&\langle \alpha_1^{(1)}+\sum_{i=2}^{24}c_i\alpha_2^{(i)}\mbox{ such that } (c_2,\ldots c_{24})\\&\mbox{ is a cyclic permutation of }(00000101001100110101111)\rangle\end{array}$ \\ \hline
$ -$    & $\Lambda_{24}\simeq L$ \\ \hline
\end{longtable}

Now let us consider a K3 surface $S$ such that $\rho(S)\geq 12$. Let us denote by $T_S$ its transcendental lattice.
We describe an algorithm which gives a $\mathcal{J}_2$-classification of the elliptic fibration on $S$.

\begin{enumerate}

\item {\bf The lattice $T$:} We define the lattice $T$ to be a
negative definite lattice such that $\rk(T)=\rk(T_S)+4$ and the
discriminant group and form of $T$ are the same as the ones of
$T_S$. The lattice $T$ is not necessarily unique. If it is not, we
choose one lattice with this property (the results obtained do not
depend on this choice). \item {\bf Assumption:} We assume that one
can choose $T$ to be a root lattice. \item {\bf The embeddings
$\phi$:} Given a Niemeier lattice $L$ we choose a set of primitive
embeddings $\phi:T\hookrightarrow L_{\rt}$ not isomorphic by an
element of the Weyl group. \item {\bf The lattices $N$ and
$N_{\rt}$: }Given a primitive embedding $\phi$ we compute the
orthogonal complement $N$ of $\phi(T)$ in $L_{\rt}$, i.e.
$N:=\phi(T)^{\perp_{L_{\rt}}}$ and $N_{\rt}$ its root lattice.
\item {\bf The lattices $W$ and $W_{\rt}$:} We denote by $W$ the
orthogonal complement of $\phi(T)$ in $L$, i.e.
$W:=\phi(T)^{\perp_{L}}$ and by $W_{\rt}$ its root lattice. We
observe that $N_{\rt}=W_{\rt}$ and $N\hookrightarrow W$ with
finite index. \item {\bf The elliptic fibration
$\mathcal{E}_{\phi}$:} The frame of any elliptic fibration on $S$
is a lattice $W$ obtained as in step 5. Moreover, the trivial
lattice of any elliptic fibration on $S$ is of the form $U\oplus
N_{\rt}=U\oplus W_{\rt}$ where $W_{\rt}$ and $N_{\rt}$ are
obtained as above. Hence, we find a $\mathcal{J}_2$-classification
of the elliptic fibration on $S$. In particular every elliptic
fibration on $S$ is uniquely associated to a primitive embedding
$\phi:T\hookrightarrow L$. Let us denote by $\mathcal{E}_{\phi}$
the elliptic fibration associated to $\phi$. \item{\bf The
singular fibers:} We already observed (cf. Section \ref{sec: basic
on elliptic fibrations}) that almost all the properties of the
singular fibers are encoded in the trivial lattice, so it is clear
that every $N_{\rt}(:=(\phi(T)^{\perp_{L_{\rt}}})_{\rt})$
determines almost all the properties of the reducible fibers of
$\mathcal{E}_{\phi}$. \item {\bf The rank of the Mordell--Weil
group:} Let $\phi$ be a given embedding. Let
$r:=\rk(MW(\mathcal{E}_{\phi}))$. Then $r=\rk
(NS(S))-2-\rk(N_{\rt})=20-\rk(T_S)-\rk(N_{\rt})$.
\item {\bf The torsion of the Mordell--Weil group:} The torsion part of the Mordell--Weil group is $\overline{W_{\rt}}/W_{\rt}(\subset W/N$) and can be computed in the following way: let $l+L_{\rt}$ be a non trivial element of $L/L_{\rt}$. If there exist $k\neq 0$ and $u \in L_{\rt}$ such that $k(l+u) \in N_{\rt}$ then $l+u \in W$ and the class of $l$ is a torsion element.
\end{enumerate}

\begin{rem}{\rm It is not always true that the lattice $T$ can be chosen to be a root lattice, and the method can be applied with some modifications without this assumption, see \cite{BKW}. Since everything is easier under this assumption and in our case we can require that $T$ is a root lattice, we described the method with the assumption (2). In particular, if $T$ is not a root lattice, then one has to consider the primitive embeddings of $T$ in $L$, but one cannot use the results in \cite[Sections 4 and 5]{Nis}, so the points (3), (4) and (5) are significantly more complicated.}
\end{rem}

\subsection{Explicit computations}\label{sec: explicit computations}
Here we apply the algorithm described in Section \ref{sec:
details} to the K3 surface $X$.
\subsubsection{Step 1}
From Proposition \ref{prop: the surface X}, we find that the transcendental lattice of $X$ is $$ T_X=\left( \begin{matrix}
 6 & 0 \\
 0 & 2
\end{matrix}
\right ).
$$
According to Nishiyama \cite{Nis}, \cite{ScSh}, \cite{BL}, $T_X(-1)$ admits a primitive embedding in $E_8$ and we can take $T$ as its orthogonal complement in $E_8$, that is
$$T=A_5 \oplus A_1.$$

\subsubsection{Step 2} We observe that $T$ is a root lattice.

\subsubsection{Step 3} We must find all the primitive embeddings  $\phi:T\hookrightarrow L_{\rt}$ not Weyl isomorphic. This has been done by Nishiyama \cite{Nis} for the primitive embeddings of $A_k$ in $A_m$, $D_n$, $E_l$ and for the primitive embeddings of $A_5\oplus A_1$ into $E_7$ and $E_8$. So we have to determine the primitive embeddings not isomorphic of $A_5\oplus A_1$ in $A_m$ and $D_n$. This will be achieved using Corollary \ref{cor: unique primitive embedding in Dn, An} and Lemma \ref{lemma: not unique in D8}. First we recall some notions used in order to prove these results.

Let $B$ be a negative-definite even lattice, let  $a\in B_{\rt}$ a root of $B$. The reflection $R_{a}$ is the isometry $R_{a}\left(  x\right) =x+\left(  a\cdot x\right)  a$  and the Weyl group of $B$, $W\left(  B\right)$, is the group generated by $R_{a}$ for $a\in B_{\rt}$.

\begin{proposition}
Let $A$ be a sublattice of $B$. Suppose there exists a sequence of
roots $x_{1},x_{2},\ldots ,x_{n}$ of $A^{\perp_B}$ with
$x_{i}\cdot x_{i+1}=\varepsilon_{i}$ and $\varepsilon_{i}^{2}=1$
then the two lattices $A\oplus x_{1}$ and $A\oplus x_{n}$ are
isometric by an element of the Weyl group of $B.$
\end{proposition}
\proof First we prove the statement for $n=2$. Since the two
sublattices $A\oplus\langle x_1\rangle$ and
$A\oplus\langle-x_1\rangle$ are isometric by $R_{x_1}$ we can
suppose that $x_1\cdot x_2=1$ (i.e. $\varepsilon_1=1$). Then
$x_1+x_2$ is also a root and is in $A^{\perp_B}.$ So the
reflection $R_{x_1+x_2}$ is equal to $I_{d}$ on $A.$ Let
$g:=R_{x_1}\circ R_{x_1+x_2}$ then $g\in W(B)$ is equal to $I_{d}$
on $A.$ Moreover $g\left( x_2\right)  =R_{x_1}\left(  x_2+\left(
\left( x_1+x_2\right) \cdot x_2\right)  \left( x_1+x_2\right)
\right) =x_1,$ and so $g$ fits. The case $n>2$ follows by
induction.\eprf

\begin{corollary}\label{cor: unique primitive embedding in Dn, An}
Suppose $n\geq9,$ $p\geq6,$ up to an element of the Weyl group
$W\left( D_{n}\right)$ or $W\left(  A_{p}\right)$ there is a
unique primitive embedding of $A_{5}\oplus A_{1}$ in $D_{n}$ or
$A_{p}.$
\end{corollary}
\proof From Nishiyama, \cite{Nis} up to an element of the Weyl
group there exists one primitive embedding of  $A_{5}$ in $D_{n}$
or $A_{p}.$ Fix this embedding.  If $M$ is the orthogonal of this
embedding then $M_{\rt}$ is $D_{n-6}$ or $A_{p-6}.$ So for two
primitive embeddings of $A_{1}$ in $M_{\rt}$ we can apply the
previous proposition.  \eprf

We study now the primitive embeddings of $A_{5}\oplus A_{1}$ in $D_{8},$ (which are not considered in the previous corollary, since the orthogonal complement of the unique primitive embedding of $A_5$ in $D_8$ is $\langle -6\rangle\oplus \langle -2\rangle^2$).

We denote by $\{\varepsilon_{i},1\leq i\leq n\}$ the canonical basis of
$\mathbb{R}^{n}.$

We can identify $D_{n}\left(  -1\right)  $ with $\mathbb{D}_{n},$ the set of
vectors of $\mathbb{Z}^{n}$ whose coordinates have an even sum.

First we recall the two following propositions, see for example
\cite{Ma}.

\begin{proposition}
The group $Aut(\mathbb{Z}^{n})$ is isomorphic to the semi-direct product
$\{\pm1\}^{n}\rtimes S_{n}$, where the group $S_{n}$ acts on $\{\pm1\}^{n}$ by
permuting the n components.
\end{proposition}

\begin{proposition}
If $n\neq4$ $\ $\ the restriction to $\mathbb{D}_{n}$ of the automorphisms of
$\mathbb{Z}^{n}$ induces an isomorphism of $Aut\left(  \mathbb{Z}^{n}\right)
$ onto the group $Aut\left(  \mathbb{D}_{n}\right)$. The Weyl group
$W\left(  \mathbb{D}_{n}\right)  $ of index two in $Aut\left( \mathbb{D}_{n}\right)  $ corresponds to those elements which induce an even number of
changes of signs of the $\varepsilon_{i}$.
\end{proposition}

\begin{lemma}\label{lemma: not unique in D8}
There are two embeddings of $A_{5}\oplus A_{1}$ in $D_{8}$ non isomorphic up to $W\left(  D_{8}\right).$
\end{lemma}
\proof Let $d_{8}=\varepsilon_{1} +\varepsilon_{2}$ and $d_{8-i+1}=-\varepsilon_{i-1}+\varepsilon_{i}$ with
$2\leq i\leq8$ a basis of $\mathbb{D}_{8}.$ We consider the embedding
\[
A_{5}\hookrightarrow\langle d_{7},d_{6},d_{5},d_{4,}d_{3}\rangle.
\]
By Nishiyama's results \cite{Nis},  this embedding is unique up to
an element of $W(D_{8})$ and we have $\left( A_{5}\right)^{\perp
_{D_{8}}}=\langle\sum_{i=1}^{6}\varepsilon
_{i}\rangle\oplus\langle x_{7}\rangle\oplus\langle d_{1}\rangle $
with $x_{7}=\varepsilon_{7}+\varepsilon_{8}.$ We see that $\pm
x_{7}$ and $\pm d_{1}$ are the only roots of $\left( A_{5}\right)
^{\perp_{D_{8}}}$.

We consider the two embeddings
\begin{align*}
A_{5}\oplus A_{1} &  \hookrightarrow\langle d_{7},d_{6},d_{5},d_{4,}d_{3}\rangle\oplus\langle x_{7}\rangle\\
A_{5}\oplus A_{1} &  \hookrightarrow\langle d_{7},d_{6},d_{5},d_{4,}d_{3}\rangle\oplus\langle d_{1}\rangle.
\end{align*}
Suppose there exists an element $R^{\prime}$ of $W(D_{8})$ such
that $R^{\prime}(x_{7})=d_{1}$ and $R^{\prime}\left(  A_{5}\right)
=A_{5},$ we shall show that $R^{\prime}\left(  d_{1}\right)  =\pm
x_{7}.$ If $z:=R^{\prime }\left(  d_{1}\right)  $ then as
$R^{\prime}$ is an isometry $z\cdot d_{1} =R'(d_1)\cdot
R'(x_7)=d_{1}\cdot x_{7}=0.$ Moreover, $z\in\left(
A_{5}\right)^{\perp _{D_{8}}}$ and so $z$ $=\pm x_{7}.$

Since $R^{\prime}\left(  A_{5}\right)  =A_{5}$, we see that $R^{\prime}|_{A_{5}}$
is an element of $O\left(  A_{5}\right)  $, the group of isometries of $A_{5}.$ We know
that $O\left(  A_{5}\right)  /W\left(  A_{5}\right)  \sim\mathbb{Z}%
/2\mathbb{Z,}$ generated by the class of $\mu:d_{7}\leftrightarrow d_{3},$
$d_{6}\leftrightarrow d_{4},$ $d_{5}\leftrightarrow d_{5}.$

Thus, we have $R^{\prime}|_{A_{5}}=\rho\in W\left(  A_{5}\right) $
or $R^{\prime}|_{A_{5}}=\rho\mu$ with $\rho\in W\left(
A_{5}\right)  .$ We can also consider $\rho$ as an element of the
group generated by reflections $R_{u}$ of $D_{8}$ with $u\in
A_{5}.$ So, for $v$ in $\left(  A_{5}\right)^{\perp_{D_{8}}}$ we
have $R_{u}(v)=v$ if $u\in A_{5}$ and then $\rho\left( v\right)
=v.$

Let $R=\rho^{-1}R^{\prime}$ then $R$ $=R^{\prime}$ on $\left(
A_{5}\right)^{\perp_{D_{8}}}.$ Since $R^{\prime}\left(
d_{1}\right) =\pm x_{7}$ and $R^{\prime}\left(  x_{7}\right)
=d_{1}$ we have $R^{\prime}|_{\langle\varepsilon
_{7},\varepsilon_{8}\rangle}=\left(
\varepsilon_{7}\rightarrow\varepsilon
_{8},\varepsilon_{8}\rightarrow-\varepsilon_{7}\right)  $ or
$\left(
\varepsilon_{7}\rightarrow\varepsilon_{7},\varepsilon_{8}\rightarrow
-\varepsilon_{8}\right)  .$ Also we have
$R|_{\langle\varepsilon_{1},\varepsilon _{2},\ldots
,\varepsilon_{6}\rangle}=I_d$ or
$\varepsilon_{i}\leftrightarrow\varepsilon _{7-i}.$

In the second case $R$ corresponds to a permutation of $\varepsilon_{i}$ with
only one sign minus; thus, $R$ is not an element of $W\left(  D_{8}\right)  .$

\subsubsection{Step 4}  For each primitive embedding of $A_5\oplus A_1$ in $L_{\rt}$, the computations of $N$ and $N_{\rt}$ are obtained in almost all the cases by \cite[Section 5]{Nis}. In the few cases not considered by Nishiyama, one can make the computation directly. The results are collected in Table \ref{table: orthogonal}, where we use the following notation.
The vectors $x_3,\, x_7,\, z_6$ in $D_n$ are defined by
\begin{eqnarray*}
x_3 & := &  d_{n-3}+ 2d_{n-2} + d_{n-1} + d_n, \\
x_7 & := &  d_{n-7}+ 2(d_{n-6} + d_{n-5} + d_{n-4} +d_{n-3} + d_{n-2})+ d_{n-1} + d_n, \\
x_7' & := &  2(d_{n-6} + d_{n-5} + d_{n-4} +d_{n-3} + d_{n-2})+ d_{n-1} + d_n, \\
z_6 & := &  d_{n-5}+2d_{n-4}+ 3d_{n-3}+ 4d_{n-2}+ 3d_{n-1} + 2d_n,\\
\widetilde{z_6} & := &  d_{n-5}+2d_{n-4}+ 3d_{n-3}+ 4d_{n-2}+ 2d_{n-1} + 3d_n,
\end{eqnarray*}
and the vectors $x,\, y$ in $E_p$  are
\begin{eqnarray*}
x  := & e_1 + e_2 + 2e_3 + 2e_4 + e_5, \
y  := & e_1 + 2e_2 + 2e_3 + 3e_4 + 2e_5 + e_6.
\end{eqnarray*}

\begin{longtable}{|c|l|l|}\caption{The orthogonal complement of the primitive embeddings $A_5\oplus A_1$ in $L_{\rt}$}\label{table: orthogonal}\\
\hline
No.  & Primitive  Embedding&  Orthogonal Complement\\
\hline
\hline
1 & $\begin{array}{r} \left\langle e_1^{(1)},\, e_3^{(1)},\ldots,\, e_6^{(1)}\right\rangle \\ \oplus\left\langle e_8^{(1)}\right\rangle\end{array}$  &
$\begin{array}{l}\left\langle \begin{array}{l} 4e_1^{(1)}+6e_2^{(1)}+8e_3^{(1)}+12e_4^{(1)}+ \\ 10e_5^{(1)}+ 8e_6^{(1)}+ 6e_7^{(1)}+3e_8^{(1)}\end{array} \right\rangle \oplus \\
 \left\langle y^{(1)}\right\rangle\oplus\left\langle e_1^{(2)},\ldots ,e_8^{(2)}\right\rangle  \oplus  \left\langle e_1^{(3)},\ldots ,e_8^{(3)}\right\rangle\end{array}$ \\
\cline{2-3}
2 & $\begin{array}{r} \left\langle e_1^{(1)},\, e_3^{(1)},\ldots,\, e_6^{(1)}\right\rangle \\ \oplus\left\langle e_1^{(2)}\right\rangle \end{array} $ &
$\begin{array}{l}\left\langle e_8^{(1)}, \begin{array}{l} 2e_1^{(1)}+3e_2^{(1)}+4e_3^{(1)}+ 6e_4^{(1)}+ \\ 5e_5^{(1)}+ 4e_6^{(1)}+  3e_7^{(1)}+2e_8^{(1)} \end{array} \right\rangle \oplus \\ \left\langle y^{(1)}\right\rangle\oplus \left\langle x^{(2)},\, e_2^{(2)},\, e_4^{(2)},\ldots,\, e_8^{(2)}\right\rangle \oplus \\ \left\langle e_1^{(3)},\ldots ,e_8^{(3)}\right\rangle \end{array}$\\
\hline
3 & $\left\langle e_1,\, e_3,\ldots,\, e_6\right\rangle\oplus\left\langle e_8\right\rangle$ &
$\begin{array}{l}\left\langle y\right\rangle\oplus\left\langle \begin{array}{l}4e_1+6e_2+8e_3+12e_4+ \\ 10e_5+8e_6+6e_7+3e_8 \end{array}\right\rangle \oplus \\  \left\langle d_1,\ldots ,d_{16}\right\rangle \end{array}$ \\
\cline{2-3}
4 & $\left\langle d_{16},\, d_{14},\ldots ,d_{11}\right\rangle \oplus \left\langle d_1\right\rangle$ &
$\left\langle z_6\right\rangle \oplus \left\langle \begin{array}{r} x_7,\, d_9,\ldots ,d_3,\\ 2d_2+d_1 \end{array}\right\rangle  \oplus  \left\langle e_1,\ldots, e_8\right\rangle$ \\
\cline{2-3}
5  & $\left\langle e_1,\, e_3,\ldots,\, e_6\right\rangle \oplus \left\langle d_{16}\right\rangle$ & $\begin{array}{l}\left\langle y\right\rangle\oplus\left\langle e_8, \, \begin{array}{r}2e_1+3e_2+4e_3+6e_4+ \\ 5e_5+4e_6+3e_7+2e_8 \end{array}\right\rangle \oplus \\
\left\langle d_{15}\right\rangle\oplus\left\langle x_3,\, d_{13},\ldots ,d_1 \right\rangle \end{array}$\\
\cline{2-3}
6 & $\left\langle d_{16},\,d_{14},\ldots ,d_{11}\right\rangle \oplus \left\langle e_1\right\rangle$ &
$\begin{array}{l}\left\langle z_6\right\rangle\oplus \left\langle x_7,\,d_9,\ldots ,d_1\right\rangle\oplus  \left\langle x,\, e_2,\, e_4,\ldots,\,e_8 \right\rangle\end{array}$\\
\hline
7 & $\begin{array}{r} \left\langle e_2^{(1)},\, e_4^{(1)},\ldots ,e_7^{(1)}\right\rangle \\ \oplus\langle e_1^{(1)}\rangle \end{array}$&
$\begin{array}{l}\left\langle \begin{array}{l}3e_1^{(1)}+4e_2^{(1)}+6e_3^{(1)}+ \\ 8e_4^{(1)}+6e_5^{(1)}+4e_6^{(1)}+2e_7^{(1)}\end{array}\right\rangle  \oplus \\ \left\langle d_1,\ldots ,d_{10}\right\rangle\oplus \left\langle e_1^{(2)},\ldots ,e_7^{(2)}\right\rangle \end{array}$ \\
\cline{2-3}
8 & $\left\langle d_{10},\,d_8,\ldots ,d_5\right\rangle \oplus \left\langle d_1\right\rangle$ & $\begin{array}{l} \left\langle z_6\right\rangle \oplus \left\langle x_7\right\rangle\oplus\left\langle d_3\right\rangle\oplus  \left\langle d_3+x_7+2d_2+d_1\right\rangle \oplus \\
\left\langle e_1^{(1)},\ldots ,e_7^{(1)}\right\rangle\oplus \left\langle e_1^{(2)},\ldots ,e_7^{(2)}\right\rangle \end{array}$\\
\cline{2-3}
9 & $\begin{array}{r}\left\langle e_1^{(1)},\, e_3^{(1)},\ldots , e_6^{(1)}\right\rangle  \\ \oplus \left\langle e_1^{(2)}\right\rangle \end{array}$ & $\begin{array}{l}\left\langle \begin{array}{l} 2e_1^{(1)}+3e_2^{(1)}+4e_3^{(1)}+ \\ 6e_4^{(1)}+5e_5^{(1)}+4e_6^{(1)}+3e_7^{(1)} \end{array} \right\rangle\oplus \\\left\langle y^{(1)}\right\rangle\oplus \left\langle x^{(2)},\, e_2^{(2)},\, e_4^{(2)},\ldots , e_7^{(2)}\right\rangle\oplus \\
 \left\langle d_1,\ldots ,d_{10}\right\rangle \end{array}$\\
\cline{2-3}
10 & $\begin{array}{r} \left\langle e_2^{(1)},\, e_4^{(1)},\ldots , e_7^{(1)}\right\rangle \\ \oplus \left\langle e_1^{(2)}\right\rangle \end{array}$ &
$\begin{array}{l}\left\langle e_1^{(1)},\, \begin{array}{l} 2e_1^{(1)}+2e_2^{(1)}+3e_3^{(1)}+ \\ 4e_4^{(1)}+3e_5^{(1)}+2e_6^{(1)}+e_7^{(1)} \end{array} \right\rangle \oplus \\ \left\langle x^{(2)},\, e_2^{(2)},\, e_4^{(2)},\ldots, e_7^{(2)}\right\rangle\oplus \left\langle d_1,\ldots ,d_{10}\right\rangle \end{array}$\\
\cline{2-3}
$11$ & $\begin{array}{r} \left\langle e_1^{(1)},\, e_3^{(1)},\ldots, e_6^{(1)}\right\rangle \\ \oplus\left\langle d_{10}\right\rangle \end{array}$ &
$\begin{array}{l}\left\langle \begin{array}{l} 2e_1^{(1)}+3e_2^{(1)}+4e_3^{(1)}+ \\ 6e_4^{(1)}+5e_5^{(1)}+4e_6^{(1)}+3e_7^{(1)} \end{array}\right\rangle \oplus \left\langle d_9\right\rangle \oplus \\ \left\langle e_1^{(2)},\ldots ,e_7^{(2)}\right\rangle\oplus \left\langle x_3,\,d_7,\ldots ,d_1\right\rangle \oplus \left\langle y^{(1)}\right\rangle\end{array}$\\
\cline{2-3}
12 & $\begin{array}{r} \left\langle e_2^{(1)},\, e_4^{(1)},\ldots , e_7^{(1)}\right\rangle\\ \oplus\left\langle d_{10}\right\rangle \end{array}$ &
$\begin{array}{l}\left\langle e_1^{(1)},\, \begin{array}{l} 2e_1^{(1)}+2e_2^{(1)}+3e_3^{(1)}+ \\ 4e_4^{(1)}+3e_5^{(1)}+2e_6^{(1)}+e_7^{(1)} \end{array} \right\rangle \oplus \\ \left\langle e_1^{(2)},\ldots ,e_7^{(2)}\right\rangle \oplus \left\langle d_9\right\rangle\oplus  \left\langle x_3,\,d_7,\ldots ,d_1\right\rangle \end{array}$\\
\cline{2-3}
13 & $\left\langle d_{10},\,d_8,\ldots ,d_5\right\rangle\oplus\left\langle e_1^{(1)}\right\rangle $ &
$\begin{array}{l}\left\langle x_7,\,d_3,\,d_2,\,d_1\right\rangle\oplus\left\langle \begin{array}{r} x^{(1)},\,e_2^{(1)},\\ e_4^{(1)},\ldots ,e_7^{(1)} \end{array}\right\rangle \oplus \\ \left\langle e_1^{(2)},\ldots ,e_7^{(2)}\right\rangle \oplus  \left\langle z_6\right\rangle\end{array}$\\
\hline
14 & $ \left\langle e_2,\,e_4,\ldots,e_6\right\rangle\oplus\langle e_1\rangle $ &
$\begin{array}{l} \left\langle \begin{array}{l} 2e_1+3e_2+4e_3+ \\ 6e_4 +5e_5+4e_6+3e_7 \end{array}\right\rangle\oplus  \left\langle a_1,\ldots ,a_{17}\right\rangle \end{array}$ \\
\cline{2-3}
15 & $\left\langle a_1,\ldots ,a_5\right\rangle\oplus \left\langle a_7\right\rangle $ & $\langle e_1, \ldots ,e_7\rangle \oplus \left\langle \begin{array}{r} a_9,\ldots ,a_{17}, \\ \sum_{j=1}^6ja_j-6a_8, \\ a_7+2a_8 \end{array}\right\rangle $ \\
\cline{2-3}
16 & $\left\langle e_1,\, e_3,\ldots, e_6\right\rangle \oplus\left\langle a_1 \right\rangle$ &
$\begin{array}{l}\left\langle y\right\rangle\oplus\left\langle \begin{array}{l} 2e_1+3e_2+4e_3+ \\ 6e_4+5e_5+4e_6+3e_7\end{array} \right\rangle\oplus \\ \left\langle a_1+2a_2,\, a_3,\ldots ,a_{17}\right\rangle \end{array}$\\
\cline{2-3}
17 & $\left\langle e_2,\, e_4,\ldots, e_7\right\rangle \oplus \left\langle a_1 \right\rangle$&
$\begin{array}{l}\left\langle \begin{array}{l} 2e_1+2e_2+3e_3+ \\ 4e_4+3e_5+2e_6+e_7 \end{array},\, e_1\right\rangle\oplus \\ \left\langle a_1+2a_2,\, a_3,\ldots ,a_{17}\right\rangle \end{array}$\\
\cline{2-3}
18 & $\left\langle a_1,\ldots,\,a_5\right\rangle\oplus\left\langle e_1\right\rangle$ & $\left\langle x,\,e_2,\,e_4,\ldots,\,e_7\right\rangle\oplus \left\langle \begin{array}{r} \sum_{j=1}^6ja_j,\\ a_7,\ldots a_{17} \end{array}\right\rangle$ \\
\hline
19 & $\left\langle d_{24},\, d_{22},\ldots ,d_{19}\right\rangle \oplus \left\langle d_1\right\rangle$ & $\left\langle z_6\right\rangle\oplus \left\langle x_7\right\rangle\oplus \left\langle x_7+d_1+2d_2,d_3,\ldots d_{17}\right\rangle$\\
\hline
$20$ & $\begin{array}{r} \left\langle d_{12}^{(1)},\,d_{10}^{(1)},\ldots ,d_7^{(1)}\right\rangle \\ \oplus\left\langle d_1^{(1)}\right\rangle \end{array}$ &
$\begin{array}{l} \left\langle \begin{array}{l} d_1^{(1)}+2d_2^{(1)}+ d_3^{(1)}+x_7^{(1)}, \\ d_3^{(1)},\, d_4^{(1)},\, d_5^{(1)} \end{array}\right\rangle  \oplus \\ \left\langle d_1^{(2)},\ldots ,d_{12}^{(2)}\right\rangle \oplus \left\langle z_6^{(1)}\right\rangle\oplus\left\langle x_7^{(1)}\right\rangle \end{array}$\\
\cline{2-3}
21 & $\begin{array}{r} \left\langle d_{12}^{(1)},\, d_{10}^{(1)},\ldots ,d_7^{(1)}\right\rangle \\ \oplus\left\langle d_{12}^{(2)}\right\rangle \end{array}$ &
$\begin{array}{r}\left\langle z_6^{(1)}\right\rangle\oplus\left\langle x_7^{(1)},\,d_5^{(1)},\ldots ,d_1^{(1)}\right\rangle\oplus \\ \left\langle d_{11}^{(2)}\right\rangle\oplus\left\langle x_3^{(2)},\, d_9^{(2)},\ldots ,d_1^{(2)}\right\rangle \end{array}$\\
\hline
22 & $\begin{array}{r} \left\langle d_7^{(1)},\,d_6^{(1)},\ldots ,d_3^{(1)}\right\rangle  \\ \oplus \left\langle d_1^{(1)}\right\rangle \end{array}$ &
$\begin{array}{l} \left\langle \widetilde{z_6}^{(1)}\right\rangle \oplus \left\langle x_7^{(1)}\right\rangle \oplus \\ \left\langle d_1^{(2)},\ldots ,d_8^{(2)}\right\rangle\oplus  \left\langle d_1^{(3)},\ldots ,d_8^{(3)}\right\rangle \end{array}$\\
\cline{2-3}
$\begin{array}{c} 22 \\ (b) \end{array}$ &
$\begin{array}{r} \left\langle d_7^{(1)},\,d_6^{(1)},\ldots ,d_3^{(1)}\right\rangle \\ \oplus \left\langle x_7^{(1)}\right\rangle \end{array}$ &
$\begin{array}{l} \left\langle \widetilde{z_6}^{(1)}\right\rangle \oplus \left\langle d_1^{(1)}\right\rangle \oplus \\ \left\langle d_1^{(2)},\ldots ,d_8^{(2)}\right\rangle\oplus \left\langle d_1^{(3)},\ldots ,d_8^{(3)}\right\rangle \end{array}$\\
\cline{2-3}
23 & $\begin{array}{r} \left\langle d_8^{(1)},\,d_6^{(1)},\ldots ,d_3^{(1)}\right\rangle  \\ \oplus \left\langle d_8^{(2)}\right\rangle \end{array}$ &
$\begin{array}{l}\left\langle z_6^{(1)}\right\rangle\oplus\left\langle x_7^{(1)}\right\rangle\oplus\left\langle d_1^{(1)}\right\rangle\oplus \left\langle d_7^{(2)}\right\rangle\oplus\\ \left\langle x_3^{(2)},\, d_5^{(2)},\ldots ,d_1^{(2)}\right\rangle \oplus\left\langle d_1^{(3)},\ldots ,d_8^{(3)}\right\rangle \end{array}$\\
\hline
24 & $\left\langle d_9,\, d_7,\ldots ,d_4\right\rangle\oplus \left\langle d_1\right\rangle$ & $\begin{array}{l} \left\langle z_6\right\rangle \oplus \left\langle x_7,\, d_1+2d_2\right\rangle  \oplus \left\langle a_1,\ldots ,a_{15} \right\rangle \end{array}$\\
\cline{2-3}
25 & $\left\langle a_1,\ldots,a_5\right\rangle\oplus\left\langle a_7\right\rangle$ &
$\begin{array}{l}\left\langle d_1,\ldots ,d_9\right\rangle\oplus  \left\langle \begin{array}{r} \sum_{j=1}^6ja_j-6a_8, \\ a_7+2a_8,\\ a_9,\, \ldots ,\,a_{15}\end{array}\right\rangle \end{array}$\\
\cline{2-3}
26 & $\left\langle d_9,\, d_7,\ldots ,d_4\right\rangle \oplus\left\langle a_1\right\rangle$ & $\left\langle z_6\right\rangle \oplus \left\langle d_1,d_2, x_7\right\rangle \oplus \left\langle  \begin{array}{r} a_1+2a_2,\\ a_3,\ldots ,a_{15} \end{array}\right\rangle$\\
\cline{2-3}
27  & $\left\langle a_1,\ldots ,a_5\right\rangle\oplus\left\langle d_9\right\rangle$ &
$\begin{array}{l} \left\langle d_8\right\rangle\oplus \left\langle x_3,\, d_6,\ldots ,d_1\right\rangle  \oplus \left\langle \begin{array}{r} \sum_{j=1}^6ja_j,\\ a_7,\ldots ,a_{15} \end{array}\right\rangle \end{array}$\\
\hline
28 & $\begin{array}{r} \left\langle e_1^{(1)},\,e_3^{(1)},\ldots,e_6^{(1)}\right\rangle \\ \oplus\left\langle e_2^{(2)}\right\rangle \end{array}$ &
$\begin{array}{l}\left\langle \begin{array}{r} e_2^{(2)}+e_3^{(2)}+2e_4^{(2)}+e_5^{(2)}, \\ e_1^{(2)},\,e_3^{(2)},\,e_5^{(2)},\, e_6^{(2)}\end{array}\right\rangle\oplus \\ \left\langle e_1^{(3)},\ldots ,e_6^{(3)}\right\rangle \oplus \left\langle e_1^{(4)},\ldots ,e_6^{(4)}\right\rangle \oplus \left\langle y^{(1)}\right\rangle \end{array}$\\
\hline
29 & $\left\langle a_1,\ldots,a_5 \right\rangle\oplus \left\langle a_7\right\rangle $ &
$\begin{array}{l}\left\langle \begin{array}{l} \sum_{j=1}^6ja_j-6a_8,\\ a_7+2a_8,\, a_9,\, a_{10},\, a_{11} \end{array} \right\rangle\oplus\\ \left\langle d_1,\ldots ,d_7\right\rangle \oplus \left\langle e_1,\ldots ,e_6\right\rangle\end{array}$ \\
\cline{2-3}
30 & $\left\langle e_1,\, e_3,\ldots, e_6\right\rangle\oplus\left\langle d_7\right\rangle$ & $\begin{array}{l} \left\langle y\right\rangle\oplus \left\langle d_6\right\rangle\oplus\left\langle x_3,\, d_4,\ldots ,d_1\right\rangle \oplus \left\langle a_1,\ldots ,a_{11}\right\rangle \end{array}$\\
\cline{2-3}
31 & $\left\langle e_1,e_3\ldots, e_6\right\rangle\oplus\left\langle a_1\right\rangle$ & $\begin{array}{l} \left\langle y\right\rangle\oplus\left\langle d_1,\ldots ,d_7\right\rangle\oplus  \left\langle a_1+2a_2,\, a_3,\ldots ,a_{11}\right\rangle \end{array}$\\
\cline{2-3}
32 & $\left\langle d_7,\,d_5,\ldots ,d_2\right\rangle\oplus \left\langle e_1\right\rangle$ & $\begin{array}{l} \left\langle z_6\right\rangle\oplus\left\langle x_7'\right\rangle \oplus \left\langle x,\, e_2,\, e_4,\,e_5,\,e_6\right\rangle  \\ \oplus \left\langle a_1,\ldots ,a_{11} \right\rangle \end{array}$\\
\cline{2-3}
33  & $\left\langle d_7,\,d_5,\ldots ,d_2\right\rangle\oplus \left\langle a_1\right\rangle$ & $\begin{array}{l} \left\langle z_6\right\rangle\oplus\left\langle x_7'\right\rangle\oplus \left\langle e_1,\ldots ,e_6\right\rangle  \oplus \\ \left\langle \begin{array}{r} a_1+2a_2,\\ a_3,\ldots ,a_{11} \end{array}\right\rangle \end{array}$\\
\cline{2-3}
34 & $\left\langle a_1,\ldots,a_5\right\rangle\oplus\left\langle d_7\right\rangle$ &
$\begin{array}{l}\left\langle d_6\right\rangle \oplus \left\langle x_3,\, d_4,\ldots,d_1\right\rangle \\ \oplus \left\langle e_1,\ldots ,e_6\right\rangle \oplus \left\langle \begin{array}{r} \sum_{j=1}^6ja_j,\\ a_7,\ldots ,a_{11}\end{array}\right\rangle \end{array}$\\
\cline{2-3}
35  & $\left\langle a_1,\ldots,a_5\right\rangle \oplus \left\langle e_1\right\rangle$ &
$\begin{array}{l} \left\langle x,\, e_2,\, e_4,\, e_5,\, e_6\right\rangle \oplus \left\langle d_1,\ldots ,d_7\right\rangle\oplus \\ \left\langle \sum_{j=1}^6ja_j, a_7,\ldots ,a_{11}\right\rangle \end{array}$\\
\hline
36 & $\begin{array}{r} \left\langle d_6^{(1)},\, d_4^{(1)},\ldots,d_1^{(1)}\right\rangle \\ \oplus\left\langle d_6^{(2)}\right\rangle \end{array}$ &
$\begin{array}{l}\left\langle z_6^{(1)}\right\rangle\oplus \left\langle d_5^{(2)}\right\rangle \oplus  \left\langle d_1^{(3)},\ldots ,d_6^{(3)}\right\rangle \oplus \\ \left\langle x_3^{(2)},\, d_3^{(2)},\,d_2^{(2)},\,d_1^{(2)}\right\rangle \oplus \left\langle d_1^{(4)},\ldots ,d_6^{(4)}\right\rangle \end{array}$ \\
\hline
37 & $\left\langle a_1^{(1)},\ldots,a_5^{(1)} \right\rangle \oplus \left\langle a_7^{(1)}\right\rangle$ &
$\begin{array}{l} \left\langle  \begin{array}{r} \sum_{j=1}^6ja_j^{(1)}-6a_8^{(1)}, \\ a_7^{(1)}+2a_8^{(1)},\, a_9^{(1)} \end{array} \right\rangle\oplus \\ \left\langle a_1^{(2)},\ldots ,a_9^{(2)}\right\rangle \oplus \left\langle d_1,\ldots ,d_6\right\rangle\end{array}$\\
\cline{2-3}
38 & $\left\langle a_1^{(1)},\ldots, a_5^{(1)}\right\rangle\oplus\left\langle a_1^{(2)}\right\rangle $ &
$\begin{array}{l}\left\langle d_1,\ldots ,d_6\right\rangle \oplus \left\langle \begin{array}{l} \sum_{j=1}^6ja_j^{(1)},\\ a_7^{(1)},\, a_8^{(1)},\, a_9^{(1)} \end{array}\right\rangle \oplus \\ \left\langle a_1^{(2)}+2a_2^{(2)},\, a_3^{(2)},\ldots ,a_9^{(2)}\right\rangle\end{array}$ \\
\cline{2-3}
39 & $\left\langle a_1^{(1)},\ldots, a_5^{(1)}\right\rangle\oplus \left\langle d_6\right\rangle$ &
$\begin{array}{l}\left\langle \sum_{j=1}^6ja_j^{(1)},\, a_7^{(1)},\,a_8^{(1)},\,a_9^{(1)}\right\rangle\oplus\left\langle d_5\right\rangle \\ \oplus \left\langle x_3,\, d_3,\,d_2,\,d_1\right\rangle  \oplus \left\langle a_1^{(2)},\ldots ,a_9^{(2)}\right\rangle\end{array}$ \\
\cline{2-3}
40 & $\left\langle d_5,\ldots ,d_1\right\rangle\oplus\left\langle a_1^{(1)}\right\rangle$ & $\left\langle z_6\right\rangle\oplus\left\langle \begin{array}{l} a_1^{(1)}+2a_2^{(1)},\\  a_3^{(1)},\ldots ,a_9^{(1)}\end{array}\right\rangle \oplus \left\langle a_1^{(2)},\ldots ,a_9^{(2)}\right\rangle$ \\
\hline
41  & $\left\langle a_1^{(1)},\ldots,a_5^{(1)}\right\rangle \oplus \left\langle a_7^{(1)}\right\rangle $ &
$\begin{array}{l}\left\langle \sum_{j=1}^6ja_j^{(1)}+3a_7^{(1)}\right\rangle \oplus \left\langle a_1^{(2)},\ldots ,a_7^{(2)}\right\rangle \\ \oplus  \left\langle d_1^{(1)},\ldots ,d_7^{(1)}\right\rangle \oplus \left\langle d_1^{(2)},\ldots ,d_7^{(2)}\right\rangle \end{array}$ \\
\cline{2-3}
42  & $\left\langle a_1^{(1)},\ldots,a_5^{(1)}\right\rangle\oplus\left\langle a_1^{(2)}\right\rangle$ &
$\begin{array}{l}\left\langle \sum_{j=1}^6ja_j^{(1)},\,a_7^{(1)}\right\rangle \oplus \left\langle \begin{array}{l} a_1^{(2)}+2a_2^{(2)},\\ a_3^{(2)},\ldots ,a_7^{(2)}\end{array}\right\rangle  \\ \oplus  \left\langle d_1^{(1)},\ldots ,d_5^{(1)}\right\rangle\oplus\left\langle d_1^{(2)},\ldots ,d_5^{(2)}\right\rangle\end{array}$\\
\cline{2-3}
43 & $\left\langle a_1^{(1)},\ldots, a_5^{(1)}\right\rangle\oplus\left\langle d_5^{(1)}\right\rangle$ &
$\begin{array}{l}\left\langle \sum_{j=1}^6ja_j^{(1)},\, a_7^{(1)}\right\rangle\oplus\left\langle d_4^{(1)}\right\rangle\oplus \\ \left\langle x_3^{(1)},\, d_2^{(1)},\,d_1^{(1)}\right\rangle  \oplus \left\langle a_1^{(2)},\ldots ,a_7^{(2)}\right\rangle\end{array}$ \\
\hline
44  & $\left\langle a_1^{(1)},\ldots,a_5^{(1)}\right\rangle\oplus\left\langle a_7^{(1)}\right\rangle$ &
$\begin{array}{l}\left\langle \sum_{j=1}^6ja_j^{(1)},\, a_7^{(1)}+2a_8^{(1)} \right\rangle\oplus \\ \left\langle a_1^{(2)},\ldots ,a_8^{(2)}\right\rangle \oplus \left\langle a_1^{(3)},\ldots ,a_8^{(3)}\right\rangle \end{array}$\\
\cline{2-3}
45 & $\left\langle a_1^{(1)},\ldots,a_5^{(1)}\right\rangle\oplus\left\langle a_1^{(2)}\right\rangle$ &
$\begin{array}{l}\left\langle \begin{array}{l} \sum_{j=1}^6ja_j^{(1)},\\ a_7^{(1)},\,  a_8^{(1)} \end{array}\right\rangle \oplus \left\langle \begin{array}{l} a_1^{(2)}+2a_2^{(2)},\\ a_3^{(2)},\ldots ,a_8^{(2)}\end{array}\right\rangle \\ \oplus \left\langle a_1^{(3)},\ldots ,a_8^{(3)}\right\rangle\end{array}$\\
\hline
46 & $\left\langle a_1,\ldots,a_5\right\rangle\oplus\left\langle a_7\right\rangle$ & $\left\langle \sum_{j=1}^6ja_j-6a_8, a_7+2a_8,\, a_9,\ldots \,a_{24}\right\rangle $ \\
\hline
47  & $\left\langle a_1^{(1)},\ldots,a_5^{(1)}\right\rangle\oplus\left\langle a_7^{(1)}\right\rangle$ &
$\begin{array}{l} \left\langle \begin{array}{r} \sum_{j=1}^6ja_j^{(1)}-6a_8^{(1)}, \\ a_7^{(1)}+2a_8^{(1)}, \\ a_9^{(1)},\ldots,a_{12}^{(1)}\end{array}\right\rangle  \oplus  \left\langle a_1^{(2)},\ldots, a_{12}^{(2)}\right\rangle \end{array}$\\
\cline{2-3}
48  & $\left\langle a_1^{(1)},\ldots,a_5^{(1)}\right\rangle\oplus\left\langle a_1^{(2)}\right\rangle$ &
$\left\langle \begin{array}{l} \sum_{j=1}^6ja_j^{(1)},\\ a_7^{(1)},\ldots ,a_{12}^{(1)} \end{array}\right\rangle\oplus \left\langle \begin{array}{l} a_1^{(2)}+2a_2^{(2)},\\ a_3^{(2)},\ldots ,a_{12}^{(2)}\end{array}\right\rangle$\\
\hline
49 & $\left\langle a_1^{(1)},\ldots, a_5^{(1)}\right\rangle\oplus\left\langle a_1^{(2)}\right\rangle$ &
$\begin{array}{l}\left\langle \begin{array}{l} a_1^{(2)}+2a_2^{(2)},\\ a_3^{(2)},\, a_4^{(2)},\, a_5^{(2)}\end{array}\right\rangle\oplus\left\langle a_1^{(3)},\ldots ,a_5^{(3)}\right\rangle \\ \oplus  \left\langle a_1^{(4)},\ldots ,a_5^{(4)}\right\rangle\oplus\left\langle d_1,\ldots ,d_4\right\rangle\end{array}$ \\
\cline{2-3}
50  & $\left\langle a_1^{(1)},\ldots, a_5^{(1)}\right\rangle\oplus\left\langle d_4\right\rangle$ &
$\begin{array}{l}\left\langle d_3\right\rangle \oplus \left\langle x_3\right\rangle \oplus\left\langle d_1\right\rangle\oplus \left\langle a_1^{(2)},\ldots ,a_5^{(2)}\right\rangle  \\ \oplus  \left\langle a_1^{(3)},\ldots ,a_5^{(3)}\right\rangle \oplus \left\langle a_1^{(4)},\ldots ,a_5^{(4)}\right\rangle\end{array}$ \\
\hline
51  & $\left\langle a_1^{(1)},\ldots, a_5^{(1)}\right\rangle\oplus\left\langle a_1^{(2)}\right\rangle$ &
$\begin{array}{l}\left\langle \sum_{j=1}^6ja_j^{(1)}\right\rangle\oplus\left\langle \begin{array}{l} a_1^{(2)}+2a_2^{(2)},\\ a_3^{(2)},\ldots ,a_6^{(2)}\end{array}\right\rangle\oplus \\ \left\langle a_1^{(3)},\ldots ,a_6^{(3)}\right\rangle \oplus \left\langle a_1^{(4)},\ldots ,a_6^{(4)}\right\rangle\end{array}$ \\
\hline
\end{longtable}

\subsubsection{Step 5 (an example: fibrations 22 and 22(b))}\label{sub sub: step 5} In order to compute $W$ we recall that $W$ is an overlattice of finite index of $N$; in fact, it contains the non trivial elements of $L/L_{\rt}$ which are orthogonal to $\phi(A_5\oplus A_1)$. Moreover, the index of the inclusion $N\hookrightarrow W$ depends on the discriminant of $N$. Indeed $|d(W)|=|d(NS(X))|=12$, so the index of the inclusion $N\hookrightarrow W$ is $\sqrt{|d(N)|/12}$.

As example we compute here the lattices $W$ for the two different
embeddings of $A_5\oplus A_1$ in $D_8$ (i.e. for the fibrations 22
and 22(b)).  Thus, we consider the Niemeier lattice $L$ such that
$L_{\rt}\simeq D_8^3$ and we denote the generators of $L/L_{\rt}$
as follows:
$v_1:=\delta_8^{(1)}+\overline{\delta_8}^{(2)}+\overline{\delta_8}^{(3)}$,
$v_2:=\overline{\delta_8}^{(1)}+\delta_8^{(2)}+\overline{\delta_8}^{(3)}$,
$v_3:=\overline{\delta_8}^{(1)}+\overline{\delta_8}^{(2)}+\delta_8^{(3)}.$

\textbf{Fibration $\# 22:$} we consider the embedding  $\varphi_1:A_5\oplus A_1\hookrightarrow L$ such that $\varphi_1(A_5\oplus A_1)=\langle d_7^{(1)},d_6^{(1)},d_5^{(1)},d_4^{(1)},d_3^{(1)}\rangle\oplus \langle d_1^{(1)}\rangle$.
The generators of the lattice $N$ are described in Table \ref{table: orthogonal} and one can directly check that $N\simeq \langle -6\rangle\oplus A_1\oplus  D_8\oplus D_8$. So, $|d(N)|=6\cdot 2^5$ and the index of the inclusion $N\hookrightarrow W$ is $2^2=\sqrt{6\cdot 2^5/12}$. This implies that there is a copy of $(\Z/2\Z)^2\subset (\Z/2\Z)^3$ which is also contained in $W$ and so in particular is orthogonal to $\varphi_1(A_5\oplus A_1)$.

We observe that $v_1$ is orthogonal to the embedded copy of
$A_5\oplus A_1$, $v_2$ and $v_3$ are not.  Moreover $v_2-v_3$ is
orthogonal to the embedded copy of $A_5\oplus A_1$. Hence $v_1$
and $v_2-v_3$ generates $W/N\simeq (\Z/2\Z)^2$. We just observe
that  $v_2-v_3\in W$ is equivalent mod $W_{\rt}$ to the vector
$w_2:=\widetilde{\delta_8^{(2)}}+\widetilde{\delta_8^{(3)}}\in W$,
so $W/N\simeq (\Z/2\Z)^2\simeq \langle v_1, w_2\rangle$. We will
reconsider this fibration in Section \ref{sec: again on 22 and
22b} comparing it
with the fibration $\# 22b$. \\

\textbf{Fibration $\# 22(b):$} we consider the other embedding of
$A_5\oplus A_1$ in $L_{\rt}$, i.e. $\varphi_2:A_5\oplus
A_1\hookrightarrow L$ such that $\varphi_2(A_5\oplus A_1)=\langle
d_7^{(1)},d_6^{(1)},d_5^{(1)},d_4^{(1)},d_3^{(1)}\rangle\oplus
\langle x_7^{(1)}\rangle.$

The generators of the lattice $N$ is described in Table
\ref{table: orthogonal} and one can directly check that $N\simeq
\langle -6\rangle\oplus A_1\oplus  D_8\oplus D_8$. As above this
implies that $W/N\simeq (\Z/2\Z)^2$ which is generated by elements
in $L/L_{\rt}$ which are orthogonal to $\varphi_2(A_5\oplus A_1)$.
In particular, $v_1-v_2$ and $v_2-v_3$ are orthogonal to
$\varphi_2(A_5\oplus A_1)$ so $v_1-v_2\in W$ and $v_2-v_3\in W$.
Moreover,
$v_2-v_3=\widetilde{\delta_8^{(2)}}+\widetilde{\delta_8^{(3)}}\mod
W_{\rt}$. So, denoted by
$w_2:=\widetilde{\delta_8^{(2)}}+\widetilde{\delta_8^{(3)}}$, we
have that $W/N\simeq \langle v_1-v_2,w_2\rangle$. We will
reconsider this fibration in Section \ref{sec: again on 22 and
22b} comparing it with the fibration $\# 22$.

\subsubsection{Step 6} We recalled in Section \ref{sec: basic on elliptic fibrations} that
each elliptic fibration is associated to a certain decomposition of the N\'eron--Severi group as a direct sum of $U$ and a lattice, called $W$.
In step 5 we computed all the admissible lattices $W$, so we classify
the elliptic fibrations on $X$. We denote all the elliptic
fibrations according to their associated embeddings; this gives
the first five columns of the Table \ref{table: main result}.
\subsubsection{Step 7}
Moreover, again in Section \ref{sec: basic on elliptic fibrations}, we recalled that each reducible fiber of an elliptic fibration is uniquely associated to a Dynkin diagram and that a Dynkin diagram is associated to at most two reducible fibers of the fibration. This completes step 7.

\subsubsection{Step 8}  In order to compute the rank of the Mordell--Weil group it suffices to perform the suggested computation, so $r=18-\rk(N_{\rt})$. This gives the sixth column of Table \ref{table: main result}.

For example, in cases 22 and 22(b), the lattice $N_{\rt}$ coincides and has rank 17, thus $r=1$ in both the cases.

\subsubsection{Step 9} In order to compute the torsion part of the Mordell--Weil group one has to identify the vectors $v\in W/N$ such that $kv\in N_{\rt}$ for a certain nontrivial integer number $k\in\Z$; this gives the last column of Table \ref{table: main result}. We will demonstrate this procedure in some examples below (on fibrations $\# 22$ and $\# 22(b)$), but first we remark that in several cases it is possible to use an alternative method either in order to completely determine $MW(\mathcal{E})_{\tors}$ or at least to bound it. We already presented the theoretical aspect of these techniques in Section \ref{subset: torsion MW theory}.

Probably the easiest case is the one where $r=0$. In this case
$MW(\mathcal{E})=MW(\mathcal{E})_{tors}$. Since $r=0$, this
implies that $\rk(N_{\rt})=18=\rk N$, so $N=N_{\rt}$. Hence
$W/N=W/N_{\rt}$, thus every element $w\in W/N$ is such that a
multiple is contained in $N_{\rt}$, i.e. every element of $W/N$
contributes to the torsion. Thus, $MW(\mathcal{E})=W/N=W/N_{\rt}$.
This immediately allows to compute the torsion for the 7 extremal
fibrations $\#2,5,10,12,28, 30, 50.$

\textbf{Fibration $\# 50\  N_{\rt}\simeq A_1^{\oplus 3}\oplus A_5^{\oplus 3}$ ($r=0$).} The lattice $N=N_{\rt}$ is $A_1^{\oplus 3}\oplus A_5^{\oplus 3}$, then $|d(N)|=2^36^3$ and $|W/N|=2\times 6$. Moreover $W/N\subset L/L_{\rt}\simeq (\Z/6\Z)^2\times \Z/2\Z$. This immediately implies that $W/N=\Z/6\Z\times\Z/2\Z$.

\textbf{Fibration $\# 1\  N_{\rt}\simeq A_1\oplus E_8^{\oplus 2}$ (fibers of special type, Proposition \ref{prop: torsion and reducible fibers}).} The presence of the lattice $E_8$ as summand of $N_{\rt}$ implies that the fibration has a fiber of type $II^*$ (two in this specific case). Hence $MW(\mathcal{E})_{\tors}$ is trivial.

\textbf{Fibration $\# 29\  N_{\rt}\simeq A_3\oplus D_7\oplus E_6$
(fibers of special type, Proposition \ref{prop: torsion and
reducible fibers}).} By Proposition \ref{prop: torsion and
reducible fibers} if a fibration has a fiber of type $IV^*$, then
the Mordell--Weil group is a subgroup of $\Z/3\Z$. On the other
hand, a fiber of type $D_7$, i.e., $I_3^*$ can only occur in
fibrations with 4 or 2-torsion or trivial torsion group. Therefore
$MW(\mathcal{E})_{\tors}$ is trivial.

\textbf{Fibration $\# 25\ N_{\rt}\simeq A_7\oplus D_9$ (the height formula, Section \ref{subset: height}).} Suppose there is a non-trivial torsion section $P$. Then, taking into account the possible contributions of the reducible fibers to the height pairing, there is $0\leq i \leq 7$ such that one of the following holds:
\[ 4= \frac{i(8-i)}{8} +1\mbox{ or }4= \frac{i(8-i)}{8} +1 + 5/4.\] After a simple calculation, one sees that neither of the above can happen and therefore the torsion group $MW(\mathcal{E})_{\tors}$ is trivial.

\textbf{Fibration $\# 22\  N_{\rt}\simeq A_1\oplus D_8\oplus D_8$.} We already computed the generators of $W/N$ in Section \ref{sub sub: step 5}, $W/N\simeq (\Z/2\Z)^2\simeq \langle v_1, w_2\rangle$.
A basis of $N_{\rt}$ is $\langle x_7^{(1)}\rangle\oplus \langle d_{i}^{(j)}\rangle_{i=1,\ldots 8, j=2,3}$. So
\begin{align*}
2v_1 & =d_1^{(1)}+2d_2^{(1)}+3d_3^{(1)}+4d_4^{(1)}+5d_5^{(1)}+6d_6^{(1)}+2d_7^{(1)}+3d_8^{(1)}\\
 &+\sum_{i=2}^3\left(d_7^{(i)}+d_8^{(i)}+2\left(\sum_{j=1}^7 d_j^{(i)}\right)\right)
\end{align*}
and $2v_1\not\in N_{\rt}$ since $d_1^{(1)}+2d_2^{(1)}+3d_3^{(1)}+4d_4^{(1)}+5d_5^{(1)}+6d_6^{(1)}+2d_7^{(1)}+3d_8^{(1)}$ is not a multiple of $x_7^{(1)}$.
Viceversa
$$2w_2\subset D_8^{(2)}\oplus D_8^{(3)}\in N_{\rt}.$$
Thus $MW(\mathcal{E})=\Z\times \Z/2\Z$.

\textbf{Fibration $\# 22(b)\  N_{\rt}\simeq A_1\oplus D_8\oplus D_8$.} Similarly, we consider the generators of $W/N\simeq (\Z/2\Z)^2\simeq \langle v_1-v_2, w_2\rangle$ computed in in Section \ref{sub sub: step 5}.
A basis of $N_{\rt}$ is $\langle d_1^{(1)}\rangle\oplus \langle d_{i}^{(j)}\rangle_{i=1,\ldots 8, j=2,3}$. So $$2w_1\not\in N_{\rt} \mbox{ and } 2w_2\in N_{\rt}.$$
Thus also in this case $MW(\mathcal{E})=\Z\times \Z/2\Z$.

\subsection{Again on fibrations $\# 22$ and $\#22(b)$}\label{sec: again on 22 and 22b}

As we can check in Table \ref{table: main result} and we proved in
the previous sections, the fibrations $\# 22$ and $\#22(b)$ are
associated to the same lattice $N$ and to the same Mordell--Weil
group. However, we proved in Lemma \ref{lemma: not unique in D8}
that they are associated to different (up to Weyl group)
embeddings in the Niemeier lattices, so they correspond to
fibrations which are not identified by the
$\mathcal{J}_2$-fibration and in particular they can not have the
same frame. The following question is now natural: what is the
difference between these two fibrations? The answer is that the
section of infinite order, which generates the free part of the
Mordell--Weil group of these two fibrations, has different
intersection properties, as we show now in two different ways and
contexts.

\textbf{Fibration  $\# 22$:} we use the notation of Section \ref{sub sub: step 5}. Moreover we fix the following notation: $\Theta_1^1:=x_7^{(1)}$ and $\Theta_i^{(j)}:=d_i^{(j)}$, $i=1,\ldots 8$, $j=2,3$ are respectively the non trivial components of the fibers of type $I_2$, $I_4^*$ and $I_4^*$ respectively.

The class $P:=2F+O-v_1$ is the class of a section of infinite
order of the fibration, generating the free part of
$MW(\mathcal{E})$ and the class $Q:=2F+O-w_2$ is the class of the
2-torsion section of the fibration. The section $P$ meets the
components $\Theta_1^1$, $\Theta_1^2$, $\Theta_1^3$ and $Q$ meets
the components $\Theta_0^1$, $\Theta_7^2$, $\Theta_7^3$. We
observe that $h(P)=3/2$ and $h(Q)=0$ which agree with
\cite[Formula 22]{ScSh} and the fact that $Q$ is a torsion section
respectively. We also give an explicit equation of this fibration
and of its sections, see \eqref{eq:22}.

\textbf{Fibration  $\# 22(b)$:} we use the notation of Section \ref{sub sub: step 5}. Moreover we fix the following notation: $\Theta_1^1:=d_1^{(1)}$ and $\Theta_i^{(j)}:=d_i^{(j)}$, $i=1,\ldots 8$, $j=2,3$ are respectively the non trivial components of the fibers of type $I_2$, $I_4^*$ and $I_4^*$ respectively.
The class $Q:=2F+O-w_2$ is the class of the 2-torsion section of the fibration. Observe that $Q$ meets the components $\Theta_0^1$, $\Theta_7^2$, $\Theta_7^3$.
The class
$$P=2F+O+v_1-v_2-\Theta_1^2-\Theta_2^2-\Theta_3^2-\Theta_4^2-\Theta_5^2-\Theta_6^2-\Theta_7^2$$
is the class of a section of infinite order, which intersects the
following components of the reducible fibers: $\Theta_1^1$,
$\Theta_7^2$, $\Theta_0^3$. This agrees with the height formula.
We also give an explicit equation of this fibration and of its
sections, see \eqref{eq:22b}.

\begin{rem}{\rm The generators of the free part of the Mordell--Weil group is clearly defined up to the sum by a torsion section. The section $P\oplus Q$ intersects the reducible fibers in the following components $\Theta_1^1$, $\Theta_0^2$, $\Theta_7^3$ (this follows by the group law on the fibers of type $I_2$ (or $III$) and $I_4^*$).}\end{rem}

\begin{rem}{\rm Comparing the sections of infinite order of the fibration 22 and the one of the fibration 22(b), one immediately checks that their intersection properties are not the same, so the frames of the elliptic fibration 22 and of elliptic fibration 22(b) are not the same and hence these two elliptic fibrations are in fact different under the $\mathcal{J}_2$-classification.

We observe that both the fibrations $\#22$ and $\#22(b)$ specialize the same fibration, which is given in \cite[Section 8.1, Table $r=19$, case 11)]{CG}. Indeed the torsion part of the Mordell--Weil group, which is already present in the more general fibrations analyzed in \cite{CG}, are the same and the difference between the fibration 22 and the fibration 22(b) is in the free part of the Mordell Weil group, so the difference between these two fibrations involve exactly the classes that correspond to our specialization.}\end{rem}

Here we also give an equation for each of the two different
fibrations \#22 and \#22(b). Both these equations are obtained from
the equation of the elliptic fibration \#8 \eqref{fib 8}. So
first we deduce an equation for \#8: Let
$c:=\frac{v}{\left(w-1\right)  ^{2}}$. Substituting $v$ by
$c\left( w-1\right)^{2}$ in \eqref{E3}, we obtain the equation of
an elliptic curve depending on $c$, which corresponds to the
fibration \#8 and with the following Weierstrass equation
\begin{equation}\label{fib 8}
E_{c}:\beta^{2}=\alpha\left(  \alpha^{2}+6c^{2}\alpha-c^{3}\left(  c-4\right)  \left(
4c-1\right)  \right).
\end{equation}

{\bf Fibration \#22} Putting $n^{\prime}=\frac{\alpha}{c^{2}\left(
4c-1\right)  },$ $\beta=\frac{yc^{2}\left( 4c-1\right)
}{4n^{\prime3}},\quad c=\frac{x}{4n^{\prime3}},$ in \eqref{fib 8}
we obtain \begin{eqnarray}\label{eq:22}
E_{n^{\prime}}:y^{2}=x\left( x^{2}-n^{\prime}\left(
n^{\prime2}-6n^(\prime)+1\right) x+16n^{\prime4}\right)
\end{eqnarray}
with singular fibers of type
$2I_{4}^{\ast}(n=0,\infty)+I_{2}+2I_{1}.$ We notice the point
$P=\left(  \left(  (n^{\prime}-1)^{2}n^{\prime},-2n^{\prime
2}(n^{\prime2}-1)\right)  \right)  $ of height $\frac{3}{2}$,
therefore $P$ and $Q=\left(  0,0\right)  $ generate the
Mordell-Weil group of $E_{n^{\prime}}$. \ To study the singular
fiber at $n^{\prime}=\infty$ we do the transformation
$N^{\prime}=\frac{1}{n},y=\frac{\beta_{1}}{N^{\prime6}},x=\frac{\alpha_{1}}{N^{\prime
4}}$ and $P=(\alpha_1,\beta_1)$ with $  \alpha_{1}=\left(  N^{\prime}-1\right)  ^{2}N^{\prime}$
and $\beta_{1}=-2N^{\prime2}\left(  N^{\prime2}-1\right)    .$ We deduce that
the section $P$ intersects the component of singular fibers at $0$ and
$\infty$ with the same subscript, so this fibration corresponds to fibration \#22.

{\bf Fibration \#22(b)}
Putting $n=\frac{2\alpha}{c\left(  4c-1\right)  },\beta=\frac{yc\left(  4c-1\right)  }%
{4n},c=\frac{-x}{2n}$ in \eqref{fib 8} we obtain
\begin{eqnarray}\label{eq:22b}
E_{n}:y^{2}=x\left(  x^{2}+2n\left(  n^{2}+3n+4\right)
+n^{4}\right)
\end{eqnarray}
with singular fibers of type
$2I_{4}^{\ast}(n=0,\infty)+I_{2}+2I_{1}.$ We notice the point
$P=\left(  4,2(n+2)^{2}\right)  $ of height $\frac{3}{2}$,
therefore $P$ and $Q=\left(  0,0\right)  $ generate the
Mordell-Weil group of $E_{n}$. Since $P$ does not meet the node of
the Weierstrass model at $n=0$, the section $P$ intersects the
component $\Theta_{0}$ of the singular fiber for $n=0,$ so this
fibration corresponds to \#22(b).

\begin{rem}{\rm Let us denote by $\mathcal{E}_9$ and $\mathcal{E}_{21}$ the elliptic fibrations \#9 and \#21 respectively.
They satisfy $Tr(\mathcal{E}_9)\simeq Tr(\mathcal{E}_{21})$ and
$MW(\mathcal{E}_9)\simeq MW(\mathcal{E}_{21})$, but
$\mathcal{E}_9$ is not $\mathcal{J}_2$-equivalent to
$\mathcal{E}_{21}$ since, as above, the infinite order sections of
these two fibrations have different intersection properties with
the singular fibers. Indeed these two fibrations correspond to
different fibrations  \cite[Case 10a) and case 10b), Section 8.1,
Table $r=19$]{CG} on the more general family of K3 surfaces
considered in \cite{CG}.}
\end{rem}


\bibliographystyle{amsplain}

\end{document}